# NONMONOTONICITY OF PHASE TRANSITIONS IN A LOSS NETWORK WITH CONTROLS

By Brad Luen, Kavita Ramanan[1] and Ilze Ziedins

*University of California, Berkeley, Carnegie Mellon University and University of Auckland*

We consider a symmetric tree loss network that supports single-link (unicast) and multi-link (multicast) calls to nearest neighbors and has capacity $C$ on each link. The network operates a control so that the number of multicast calls centered at any node cannot exceed $C_V$ and the number of unicast calls at a link cannot exceed $C_E$, where $C_E, C_V \leq C$. We show that uniqueness of Gibbs measures on the infinite tree is equivalent to the convergence of certain recursions of a related map. For the case $C_V = 1$ and $C_E = C$, we precisely characterize the phase transition surface and show that the phase transition is always nonmonotone in the arrival rate of the multicast calls. This model is an example of a system with hard constraints that has weights attached to both the edges and nodes of the network and can be viewed as a generalization of the hard core model that arises in statistical mechanics and combinatorics. Some of the results obtained also hold for more general models than just the loss network. The proofs rely on a combination of techniques from probability theory and dynamical systems.

**1. Introduction.** In Section 1.1 we describe and motivate the basic loss network model analyzed in this paper. In Section 1.2 we summarize the main results and place our results in the context of related prior work. An outline of the paper is provided in Section 1.3.

1.1. *Description of the model and motivation.* Many stochastic processes that arise in applications take values in $S^V$, where $S$ is a finite state space

Received July 2005; revised February 2006.
[1]Supported in part by the NSF Grants DMS-04-06191, DMI-03-23668-0000000965 and DMS-04-05343.

and $V$ is the vertex set of a finite or countably infinite loopless graph $G = (V, E)$ with edge set $E$. A generic example is the Ising model of ferromagnetism [1], where $S = \{-1, 1\}$ represents the set of possible magnetic spins at a vertex site of $G$, and configurations in $S^V$ that have opposing spins on neighboring sites are discouraged (i.e., assigned lower probability). Another well-known example is the hard core model, which was first proposed to study the equilibrium behavior of a lattice gas consisting of particles with nonnegligible radii that cannot overlap [1, 14, 23]. Here, the state $s \in S = \{0, 1\}$ of a node represents the occupation number of a particle, and particles cannot occupy neighboring vertices on the lattice. In other words, the set of feasible configurations $\Omega_G^{\mathrm{hc}}$ is given by

$$\Omega_G^{\mathrm{hc}} \doteq \{\sigma \in S^V : \sigma_x + \sigma_y \leq 1 \ \forall \, x, y \in V \text{ such that } xy \in E\}.$$

The hard core model is an example of a stochastic process subject to "hard constraints," in which certain configurations are forbidden (rather than just assigned lower probability, as in the case of the Ising model). Processes with such "hard constraints" arise in fields as diverse as statistical mechanics, combinatorics and telecommunications. In particular, the hard core model also arises in the study of random independent sets of a graph [5, 12, 20] and in the analysis of multicasting on telecommunication loss networks [18, 26]. The consideration of different applications has led to the study of various generalizations of the hard core model [2, 5, 13, 22, 24, 26, 28, 29]. One natural generalization of the hard core model arises from the study of stochastic loss networks in telecommunications.

A general loss network consists of $E$ links (or resources), indexed by $e$, with the $e$th link having integer capacity $C_e$. Calls are indexed by a set $\mathcal{R}$, and type $r \in \mathcal{R}$ calls arrive as a Poisson process with rate $\gamma_r$ and stay in the network for a random duration. Upon arrival, a call of type $r$ requests capacity $A_{er}$ from the $e$th link for each $e$, where the set of links $\{e : A_{er} > 0\}$ is often referred to as the *route* of the call. If the requested capacity is available simultaneously on all links along its route, then the call is accepted; otherwise it is rejected and lost. If the call is accepted, it occupies the requested capacity for its duration. The stationary distribution for a loss network depends on the call duration distribution only through its mean [7], so we allow call durations to have a general distribution with finite mean, which, without loss of generality, can be assumed to be equal to 1. All arrival processes and call durations are assumed to be independent of one another. For an excellent introduction to loss networks, see the review paper by Kelly [18].

In this work, we consider loss networks with the structure of a graph $G = (V, E)$ that is a regular tree (though the description in this section holds for a general graph $G$). All links or edges in the network have the same capacity $C$, and the network supports two classes of calls—multicast calls that arrive at a node with a Poisson rate $\nu$ and require one unit of capacity from each of the edges incident to it, and unicast calls that arrive at an edge with a Poisson rate $\lambda$ and occupy one unit of capacity on that edge. In addition to the capacity constraints that determine whether or not a call can be accepted, we also impose a simple control, whereby a call is accepted only if, after acceptance, there will be no more than $C_V \leq C$ multicast calls at any node and no more than $C_E \leq C$ unicast calls at any link in the network. Note that when $C_V = C_E = C$, the additional admission control is redundant and the controlled unicast-multicast model reduces to the uncontrolled unicast-multicast model that was introduced in [26].

Let $\mathbf{n}(s) \doteq \{n_v(s), n_e(s), v \in V, e \in E\}$, where $n_v(s)$ and $n_e(s)$ represent the numbers of multicast and unicast calls in progress at time $s$ at node $v$ and on edge $e$, respectively. Then, given any graph $G = (V, E)$, under the above admission rule for every $s \in (0, \infty)$, $\mathbf{n}(s)$ takes values in the set $\hat{\Omega}_G$ of feasible configurations given by

$$(1.1) \quad \hat{\Omega}_G \doteq \{\omega \in S_{C_V}^V \times S_{C_E}^E : \omega_u + \omega_{uv} + \omega_v \leq C \ \forall u, v \in V \text{ such that } uv \in E\},$$

where for any positive integer $C$,

$$(1.2) \quad S_C \doteq \{0, 1, \ldots, C\}.$$

Given vectors $\vec{\nu} = (\nu_0, \nu_1, \ldots, \nu_{C_V}) \in \mathbb{R}_+^{C_V+1}$ and $\vec{\lambda} = (\lambda_0, \lambda_1, \ldots, \lambda_{C_e}) \in \mathbb{R}_+^{C_E+1}$, we define $\hat{m}_G^{\vec{\nu}, \vec{\lambda}}$ to be the following probability distribution on $\hat{\Omega}_G$:

$$(1.3) \quad \hat{m}_G^{\vec{\nu}, \vec{\lambda}}(\omega) \doteq \begin{cases} \dfrac{\prod_{v \in V} \nu_{\omega_v} \prod_{e \in E} \lambda_{\omega_e}}{Z_G^{\vec{\nu}, \vec{\lambda}}}, & \text{if } \omega \in \hat{\Omega}_G, \\ 0, & \text{otherwise}, \end{cases}$$

where $Z_G^{\vec{\nu}, \vec{\lambda}}$ is the usual normalizing constant defined by

$$Z_G^{\vec{\nu}, \vec{\lambda}} = \sum_{\omega \in \hat{\Omega}_G} \prod_{v \in V} \nu_{\omega_v} \prod_{e \in E} \lambda_{\omega_e}.$$

Note that the quantity $\prod_{v \in V} \nu_{\omega_v} \prod_{e \in E} \lambda_{\omega_e}$ is proportional to the probability of the configuration $\omega$. We will refer to $\vec{\nu}$ and $\vec{\lambda}$ as *weight vectors* (they are sometimes also called *activity vectors* in the literature).

For finite $G$, it is well known [18] that the process $\mathbf{n}(\cdot)$ has a unique stationary distribution, which is supported on $\hat{\Omega}_G$ and is of product-form, given explictly by $\hat{m}_G^{\vec{\nu}, \vec{\lambda}}$, with weight vectors

$$(1.4) \quad \begin{aligned} \nu_i &\doteq \frac{\nu^i}{i!} \qquad \text{for } i = 0, 1, \ldots, C_V, \\ \lambda_i &\doteq \frac{\lambda^i}{i!} \qquad \text{for } i = 0, 1, \ldots, C_E. \end{aligned}$$

Here, we adopt the usual convention that $0! = 1$, so that $\nu_0 = \lambda_0 = 1$. The structure of the weight vectors given in (1.4) is a consequence of the assumption that calls arrive as a Poisson process and are blocked if they cannot be served immediately (there is no queueing), with each accepted call allocated one unit of capacity. Thus, the downward rate of calls at a node with occupancy $i$ is exactly $i$ (since call durations have mean 1). On the other hand, networks of interest often support processor sharing rather than allocate a fixed capacity to each accepted call. In the latter case, an arriving call is accepted only if the buffer capacity $C$ is not exceeded at every link along which the call requires capacity (so that the set of feasible configurations remains the same), but all multicast calls present at a node are served in a processor-sharing manner by a processor sitting at the node and, likewise, all unicast calls present at a link are served in a processor-sharing manner by a processor at the link. Motivated by the processor-sharing discipline, in the treatment below, we will also consider product-form distributions of the form (1.3) with the weight vectors $\vec{\nu}$ and $\vec{\lambda}$ defined instead by

$$(1.5) \quad \nu_i \doteq \nu^i \quad \text{for } i = 0, 1, \ldots, C_V \quad \text{and} \quad \lambda_i \doteq \lambda^i \quad \text{for } i = 0, 1, \ldots, C_E.$$

Such weights have also been considered for the uncontrolled pure multicast model on the lattice in [24]. Another purpose of considering alternative weights is to gain some insight into the influence of different weight vectors on the existence of phase transitions.

More generally, we will refer to any model possessing given weight vectors $\vec{\nu}$, $\vec{\lambda}$ and controls $C_E$, $C_V$ as the *controlled unicast-multicast model*, recognizing that, in the special cases where the weights are given by (1.4) and (1.5), these have a specific interpretation in the communications context as models for loss and processor sharing models, respectively.

Although the product measure (1.3) has an explicit representation, for large graphs, the exact computation of the stationary probabilities is made difficult by the combinatorial complexity of calculating the normalizing constant. Instead, in such cases, one often studies the behavior of the measures in the asymptotic limit as the size of the graph increases. The appropriate measures on unicast-multicast configurations on the infinite limit graph $G = (V, E)$ that provide insight into the behavior of product-form measures on sufficiently large finite graphs are Gibbs measures. A Gibbs measure is a distribution of a countably infinite family of random variables which admits some prescribed conditional probabilities [14]. Roughly speaking, the Gibbs measures studied in this work are characterized by the property that the distribution of the configuration on any finite subset $U$ of $V \cup E$, conditioned on the configuration on the complement, is equal to the regular conditional probability of an associated product-form measure of the configuration on $U$, given the configuration on the boundary of $U$ (see Definition 3.2 for a more precise formulation). Unlike stationary distributions on finite graphs, the associated Gibbs measures on infinite graphs need not be unique. In particular, there may be multiple Gibbs measures associated with a given pair of weight vectors. A *phase transition* is said to occur if one can identify the boundaries between weights for which there is a unique Gibbs measure and weights for which there are multiple Gibbs measures. Moreover, the phase transition is said to be *monotone* with respect to a certain parameter if the existence of multiple Gibbs measures at a certain value $s$ of the parameter implies the existence of multiple Gibbs measures for all parameter values greater than $s$. (Note that, since one can always consider the reciprocal of the parameter, there is no loss of generality in assuming that multiple Gibbs measures persist only with an *increase* in the parameter.) An illuminating discussion of the implications of phase transitions on the infinite graph for the nature of the stationary distribution on large finite graphs can be found in [19].

An important motivation for studying this model is to gain a better understanding of the general structure of Gibbs measures associated with product-form distributions on graphs that have weights on edges in addition to nodes. Such measures have received much less attention than those where weights are just on the nodes, which have been studied, for example, in [5, 13, 17, 28, 29]. The uncontrolled model in the presence of unicast and multicast calls (with $C_V = C_E = C$) for the simple case $C = 1$ was studied in [26]. As discussed in further detail in Section 6, we expect that the controlled model here (with $C_V = 1$) may provide insight into the behavior of the uncontrolled unicast-multicast model when $C > 1$.

Interest in this model initially arose from the need to gain insight into the effects of multicasting on the performance of a large network. Multicasting is a major recent innovation in communication networks, which was introduced to deal with the emergence of new multimedia applications that require simultaneous connections between multiple users (as in audio or video conferencing) and simultaneous transmissions to multiple users (as in the Internet) [8, 25]. In [22, 26] we showed that models with multicasting can give rise to phase transitions. The occurence of such phase transitions is undesirable, because customers experience this as unfairness in the network (in the form of spatially heterogeneous blocking probabilities, even though arrival rates are homogeneous). It is natural to consider whether controls can be used to mitigate this effect, and the control we propose in this paper is a particularly simple one. It has a similar flavor to traditional trunk reservation, but because reversibility is preserved, the stationary distribution still has product-form (unlike the case of trunk reservation).

1.2. *Main results and related prior work.* The first main result of the paper, Theorem 1.1, shows that, given a pair of weight vectors $\vec{\lambda}$ and $\vec{\nu}$, uniqueness of Gibbs measures for the controlled unicast-multicast model described above is equivalent to the convergence of certain recursions of the associated map $\Phi^{\vec{\nu},\vec{\lambda}}$ introduced below. Given positive integers $C_V$, $C_E$ and $C$ and weight vectors $\vec{\nu}$ and $\vec{\lambda}$, the mapping $\Phi^{\vec{\nu},\vec{\lambda}} \colon \mathbb{R}_+^{C_V} \to \mathbb{R}_+^{C_V}$, where $\Phi^{\vec{\nu},\vec{\lambda}}(\xi) = (\Phi_1^{\vec{\nu},\vec{\lambda}}(\xi), \ldots, \Phi_{C_V}^{\vec{\nu},\vec{\lambda}}(\xi))$, is defined by

$$(1.6) \qquad \Phi_k^{\vec{\nu},\vec{\lambda}}(\xi) \doteq \nu_k \left( \frac{\sum_{i=0}^{\min(C-k,C_E)} \lambda_i [1 + \sum_{j=1}^{\min(C-k-i,C_V)} \xi(j)]}{\sum_{i=0}^{C_E} \lambda_i [1 + \sum_{j=1}^{\min(C-i,C_V)} \xi(j)]} \right)^q$$

for $k = 1, \ldots, C_V$ and $\xi \in \mathbb{R}_+^{C_V}$. We refer to the mapping $\Phi^{\vec{\nu},\vec{\lambda}}$ as the *random field map* associated with the weight vectors $\vec{\nu}$ and $\vec{\lambda}$. We now state Theorem 1.1. For a precise definition of Gibbs measures associated with the controlled model, refer to Definition 3.2. In the following, $\mathbf{0}$ denotes the zero vector:

THEOREM 1.1. *Let $\Phi^{\vec{\nu},\vec{\lambda}}$ be the random field map defined in (1.6). Then, the unicast-multicast model on the $(q+1)$-regular tree, $G$, has a unique Gibbs measure if and only if, given the initial condition $\xi^{(0)} = \mathbf{0} \in \mathbb{R}_+^{C_V}$, the sequence defined iteratively by $\xi^{(n)} \doteq \Phi^{\vec{\nu},\vec{\lambda}}(\xi^{(n-1)})$ converges to a limit. In this case, the limit will be the unique fixed point $\xi_*^{\vec{\nu},\vec{\lambda}}$ of the map $\Phi^{\vec{\nu},\vec{\lambda}}$.*

The proof of Theorem 1.1 is given at the end of Section 3.2. Similar correspondences have been obtained for the case of Gibbs measures on a feasible subset of the configuration space $S^V$ for some finite set $S$, when the weight vectors are attached to just the nodes $v \in V$ of a regular tree (see [27] for the case of a binary state space, [28, 29] for the general case of bounded, positive Markov specifications on a countable state space and [5] for more recent work on the case of product-form Markov specifications). However, when the edges also have weights and the set of feasible configurations is given by (1.1), then the graph with respect to which the specifications are Markov is no longer a tree and, thus, the results of [5, 28, 29] are no longer directly applicable. Theorem 1.1 extends this correspondence to the more complicated setting in the presence of unicast calls—where the edges also have weights, and the specifications are Markov with respect to the augmented graph $\hat{G}$ defined in Section 3.1.1, rather than with respect to the regular tree $G$.

Analysis of recursions for the uncontrolled unicast-multicast model was carried out in [26]. In particular, Theorem 1.1 also shows that the results of [26] have implications for the corresponding Gibbs measures. Specifically, it rigorously proves that when $C = 1$, the value obtained in (3.7) of [26] is indeed the phase transition point for the (uncontrolled) unicast-multicast model and, also, that phase transitions continue to occur for the hard core model on the tree, even when the Markov specifications are nonhomogeneous and parity-dependent (see Theorem 5.1 of [26]). The latter property is in contrast to the behavior of the hard core model on the lattice, where phase transitions disappear when the Markov specifications are nonhomogeneous [2, 15].

The second important result of this paper is a precise characterization of the phase transition surface (namely the boundary between where there is uniqueness of Gibbs measures and where there are multiple Gibbs measures) for the controlled unicast-multicast model on the $q + 1$-regular tree with $C_V = 1$, $C_E = C$ and weight vectors $\vec{\lambda}$ that satisfy Assumption 1.1 stated below. Note that, when $C_V = 1$, since $\nu_0 = 1$, the weight vector $\vec{\nu}$ can be identified with the scalar $\nu = \nu_1$.

ASSUMPTION 1.1. For every $C \geq 2$, the weight vector $\vec{\lambda} = (\lambda_0, \lambda_1, \ldots, \lambda_C)$ satisfies

$$(1.7) \qquad \Lambda_{C-1}^2 - \Lambda_C \Lambda_{C-2} > 0,$$

where $\Lambda_C \doteq \sum_{i=0}^{C} \lambda_i$.

As shown in Lemma 4.1, this assumption is satisfied, in particular, by the weight vectors of interest described in (1.4) and (1.5). The set $\mathcal{A}$, the mapping $J$ and the quadratic $Q$ defined below will be useful in characterizing the phase transition surface.

DEFINITION 1.1. We define $\mathcal{A}$ to be the subset of the $(q, C, \vec{\lambda})$-parameter space $\{2, 3, \ldots\}^2 \times (\bigcup_{d=1}^{\infty} \mathbb{R}_+^d)$ for $(q, C, \vec{\lambda})$ that is characterized by the condition that $(q, C, \vec{\lambda}) \in \mathcal{A}$ if and only if

$$(q+1)^2 \Lambda_C \Lambda_{C-2} < (q-1)^2 \Lambda_{C-1}^2. \qquad \text{(Condition A)}$$

The mapping $J : [0, \infty) \to [0, \infty)$ is defined by

$$(1.8) \qquad J(\xi) \doteq \xi \left( \frac{\Lambda_C + \xi \Lambda_{C-1}}{\Lambda_{C-1} + \xi \Lambda_{C-2}} \right)^q$$

and the quadratic polynomial $Q$ is equal to

$$(1.9) \qquad \begin{aligned} Q(\alpha) = {} &\alpha^2 \Lambda_{C-1} \Lambda_{C-2} \\ &+ \alpha[(1-q)\Lambda_{C-1}^2 + (1+q)\Lambda_C \Lambda_{C-2}] + \Lambda_C \Lambda_{C-1}. \end{aligned}$$

The result can now be summarized in the following theorem, the proof of which is given at the end of Section 4.1:

THEOREM 1.2. *Given $\vec{\lambda}$ satisfying Assumption 1.1, the controlled unicast-multicast model with $C_V = 1$ and $C_E = C$ on the $(q+1)$-regular tree associated with the parameters $(q, C, \vec{\lambda}, \nu)$ has multiple Gibbs measures if and only if $(q, C, \vec{\lambda}) \in \mathcal{A}$ and $\nu \in (\nu_-, \nu_+)$, where $\nu_- = J(\alpha_-)$ and $\nu_+ = J(\alpha_+)$, with $J$ being the mapping defined in (1.8) and $\alpha_-$ and $\alpha_+$ being the distinct positive real roots of the quadratic $Q$ specified in (1.9).*

Note that the above theorem shows that whenever the phase transition occurs, it is nonmonotone in the parameter $\nu$, for fixed $(q, C, \vec{\lambda})$. In general, determining monotonicity of Gibbs measures is a nontrivial task. Indeed, the question of monotonicity of phase transitions for even the hard core model on certain graphs such as the cubic lattice $\mathbb{Z}^n$ remains an open problem. Interesting recent results have shown that the hard core model is monotone on certain $n$-dimensional lattices [16], but nonmonotone on certain other graphs [4]. The graph on which nonmonotonicity was demonstrated in [4] is a somewhat artificially modified tree, having several pendant nodes hanging from each node. In our work, nonmonotonicity is shown to occur in a related model arising in the context of a real application.

1.3. *Outline of the paper.* The outline of the paper is as follows. In Section 2 we relate the random field map $\Phi^{\vec{\nu}, \vec{\lambda}}$ to recursions that characterize the limiting stationary distribution at a central node of a large finite tree network and obtain expressions for the blocking probabilities. In Section 3 we provide some general background on Gibbs measures, establish properties of the Gibbs measures associated with the controlled unicast-multicast model and present the proof of Theorem 1.1, which reduces the study of uniqueness of Gibbs measures to the analysis of recursions of the random field map. In Section 4 we analyze these recursions and study the phase transitions for the case $C_V = 1$ and $C_E = C$. Numerical results are presented in Section 5 and a concluding discussion is given in Section 6.

**2. Limiting stationary distributions.** In this and subsequent sections, we only consider the case when the graph $G$ has the form of a $(q+1)$-regular tree $T$ [recall that a $(q+1)$-regular tree is the unique graph with no cycles such that each node has precisely $q+1$ neighbors]. Although the computation of the normalizing constant $Z_G^{\vec{\nu}, \vec{\lambda}}$ and, therefore, of the stationary distribution, is hard for general large graphs, the calculation is somewhat simplified when the graph is a tree. In particular, due to the absence of cycles in a tree, it is possible to derive a recursion relation for the normalizing constants for trees of increasing size (see [26] for the recursion relation for the uncontrolled version of this model).

Let $T_m$ be the complete $q$-ary tree of height $m$, with the root node denoted by $O$. In other words, $T_m$ is a rooted tree in which the root node, $O$, has degree $q$, the $qm$ terminal or leaf nodes have degree 1 and all other nodes have degree $q + 1$. Fix $\vec{\nu}$ and $\vec{\lambda}$ and let $\Omega_m \doteq \hat{\Omega}_{T_m}$. In order to derive the recursion relation, we take advantage of the property that the removal of any edge in a tree splits that tree into two disjoint trees. For $i \in \{0, 1, \ldots, C_V\}$, let $Z_m(i)$ be the weighted sum of all feasible configurations on $T_m \subset T$ that have $i$ multicast calls at the root node $O$. More precisely, define

$$Z_m(i) \doteq \sum_{n \in \Omega_m \,:\, n_0 = i} \prod_{v \in V} \nu_{n_v} \prod_{e \in E} \lambda_{n_e}.$$

Then, the vector $(Z_{m+1}(i), i = 0, 1, \ldots, C_V)$ for the tree $T_{m+1}$ can be expressed in terms of that for $T_m$ as follows:

$$(2.1) \qquad Z_{m+1}(i) = \nu_i \left( \sum_{j=0}^{\min(C-i, C_E)} \lambda_j \sum_{k=0}^{\min(C-i-j, C_V)} Z_m(k) \right)^q.$$

Since a finite $(q + 1)$-regular spherical tree $T$ (say of diameter $2L$) can be decomposed into a graph containing the central node 0 with $q + 1$ edges incident to it and $q + 1$ rooted trees each of height $L - 1$, the normalizing constant for the tree $T$ is given by

$$Z_T^{\vec{\nu}, \vec{\lambda}} = \sum_{i=0}^{C_V} Z_T^{\vec{\nu}, \vec{\lambda}}(i),$$

where

$$Z_T^{\vec{\nu}, \vec{\lambda}}(i) = \nu_i \left( \sum_{j=0}^{\min(C-i, C_E)} \lambda_j \sum_{k=0}^{\min(C-i-j, C_V)} Z_{L-1}(k) \right)^{q+1}.$$

Given boundary conditions for the external nodes, the normalizing constant $Z_T^{\vec{\nu}, \vec{\lambda}}$ can be calculated from the recursion (2.1) and the expressions in the last two displays. Thus, the problem of computing the stationary distribution of a finite tree is reduced to that of analyzing the recursion relation (2.1). In Section 3 (see the proof of Theorem 1.1) we show that the recursion relations are also useful in determining the structure of Gibbs measures on the infinite $(q + 1)$-regular tree.

The components of the vector $(Z_m(i), i \in \{0, 1, \ldots, C_V\})$ tend to infinity as $m \to \infty$. However, the stationary probability that the occupancy at the root of $T_m$ is equal to $i$, where $i \in \{0, 1, \ldots, C_V\}$, is bounded and lies in $[0, 1]$. Thus, it is more convenient to consider a renormalized vector, which carries adequate information for computing the stationary distribution. For $i = 0, \ldots, C_V$ and $m \in \mathbb{N}$, define

$$\xi_m(i) \doteq \frac{Z_m(i)}{Z_m(0)}$$

and let $\xi_m \doteq (\xi_m(1), \ldots, \xi_m(C_V))$. Note that $\xi_m$ is well defined since $Z_m(0)$ is strictly positive for every $m$. The recursion (2.1) can then be recast into the more convenient form

$$\xi_{m+1} \doteq \Phi^{\vec{\nu}, \vec{\lambda}}(\xi_m),$$

where, as defined in (1.6), $\Phi^{\vec{\nu}, \vec{\lambda}}$ is the random field map associated with the weight vectors $\vec{\nu}$ and $\vec{\lambda}$.

For a loss network, important performance measures of interest are the blocking probabilities of different types of calls. Since calls arrive as Poisson processes, the stationary blocking probability of a call is equal to the stationary probability that the system is in a configuration state such that the addition of the call to the network would take the configuration out of the feasible set $\hat{\Omega}_G$. For instance, if $\beta_v^L$ is the blocking probability of a multicast call arriving to node $v$ on the tree $T_L$, we can write

$$1 - \beta_v^L = \frac{\sum_{i=0}^{C_V-1} \nu_i (\sum_{j=0}^{\min(C-i-1,C_E)} \lambda_j \sum_{k=0}^{\min(C-i-j-1,C_V)} \xi_{L-1}(k))^{q+1}}{\sum_{i=0}^{C_V} \nu_i (\sum_{j=0}^{\min(C-i,C_E)} \lambda_j \sum_{k=0}^{\min(C-i-j,C_V)} \xi_{L-1}(k))^{q+1}}.$$

An expression for the blocking probabilities of unicast calls can also be derived in a similar fashion. From these expressions, it is possible to deduce that, in contrast to multicast blocking probabilities, the unicast blocking probabilities are always spatially homogeneous, in the sense that they depend only on the distance from the center of the network and, in the limit as $L \to \infty$, the links at the center of the network have the same blocking probability.

## 3. Gibbs measures.

In this section we introduce Gibbs measures for the controlled unicast-multicast model on the infinite $(q+1)$-regular tree described in Section 1.1. Section 3.1.1 introduces some basic notation and Section 3.1.2 provides general background on Gibbs measures and Markov specifications. Section 3.2 focuses on Gibbs measures associated with the controlled unicast-multicast model, with the main results being Theorem 3.1 and the proof of Theorem 1.1.

### 3.1. Notation and definitions.

#### 3.1.1. Basic notation.
In this section, let $G = (V, E)$ be the graph representing a regular infinite tree network with nodes or vertices $V$ and links or edges $E$. We will occasionally write $v \sim w$ if $vw \in E$. Let $\hat{G} = (\hat{V}, \hat{E})$ be the augmented graph defined by

$$\hat{V} \doteq V \cup E,$$

$$\hat{E} \doteq E \cup \{xe : x \in V, e \in E \text{ and } xz = e \text{ for some } z \in V\}.$$

For any $U \subset V$, let $\partial U \doteq \{x \in V \setminus U : xz \in E \text{ for some } z \in U\}$ and $\overline{U} \doteq U \cup \partial U$ denote, respectively, the boundary and closure of $U$ with respect to the neighbor relation in $E$. Likewise, for any $U \subset \hat{V}$, let $\hat{\partial} U \doteq \{x \in \hat{V} \setminus U : xz \in \hat{E} \text{ for some } z \in U\}$ and $\Delta U \doteq U \cup \hat{\partial} U$ denote, respectively, the boundary and closure of $U$ associated with the neighbor relation in $\hat{E}$. For $U \subset \hat{V}$, for conciseness we will often use $U_V$ to denote $U \cap V$ and $U_E$ to denote $U \cap E$. We will also write $V$ (resp., $E$) for the set of nodes in $\hat{V}$ that correspond to nodes (resp., edges) in $G$. In general, if $A \subset V$, we will let $A \subset \hat{V}$ denote the set of nodes in $\hat{V}$ corresponding to the nodes in $A \subset V$, and likewise for $B \subset E$. For $U \subset V$, we use $G[U]$ to denote the graph with the set of vertices equal to $U$ and the set of edges comprising $xy \in E$, such that $x, y \in U$. For $U \subset \hat{V}$, $\hat{G}[U]$ is defined analogously.

For $U \subset V$, $\sigma_U : S_{C_V}^V \to S_{C_V}^U$ is the projection $\sigma_U(\omega) = \{\omega_i, i \in U\}$ onto $U$, and $\sigma_V$ is simply denoted by $\sigma$. Analogously for $U \subset \hat{V}$, $\hat{\sigma}_U : S_{C_V}^V \times S_{C_E}^E \to S_{C_V}^{U_V} \times S_{C_E}^{U_E}$ is the projection $\hat{\sigma}_U(\omega) = \{\omega_i, i \in U\}$ onto $U$ and, again, the identity mapping $\hat{\sigma}_{\hat{V}}$ is denoted simply by $\hat{\sigma}$. Moreover, for convenience, $\sigma_i$ and $\hat{\sigma}_i$ will be used instead of $\sigma_{\{i\}}$ and $\hat{\sigma}_{\{i\}}$, respectively. Let $\mathcal{V}$ and $\hat{\mathcal{V}}$ be the set of finite nonempty subsets of $V$ and $\hat{V}$, respectively. For $U, W \subset V$ with $U \cap W = \varnothing$, $\tau \in S_{C_V}^U$ and $\eta \in S_{C_V}^W$, $\tau \eta = \{\tau_x : x \in U\} \cup \{\eta_x : x \in W\} \in S_{C_V}^{U \cup W}$ represents the concatenation of $\tau$ and $\eta$, and the obvious analogue holds for configurations in $S_{C_V}^V \times S_{C_E}^E$. For $U \subset V$, let $\mathcal{F}(U)$ be the $\sigma$-field in $S_{C_V}^U$ generated by sets of the form $\{\sigma_i = \tau\}$ for some $i \in U$ and $\tau \in S_{C_V}$. Likewise, for $U \subset \hat{V}$, let $\hat{\mathcal{F}}(U)$ be the $\sigma$-field in $S_{C_V}^{U_V} \times S_{C_E}^{U_E}$ generated by sets of the form $\{\hat{\sigma}_i = \tau\}$, where either $i \in U_V$ and $\tau \in S_{C_V}$, or $i \in U_E$ and $\tau \in S_{C_E}$.

For conciseness, henceforth in this section we will use $S$ to denote $S_{C_V}$ and, for $U \subset \hat{V}$, we will use $\hat{S}^U$ to denote $S_{C_V}^{U_V} \times S_{C_E}^{U_E}$. The set given in (1.1) of feasible configurations for the controlled unicast-multicast model on any graph $G$ with vertex set $V$ and edge set $E$ can then be rewritten as

$$(3.1) \quad \hat{\Omega}_G \doteq \{\omega \in \hat{S}^{\hat{V}} : \omega_x + \omega_{xy} + \omega_y \leq C \text{ for all } x, y \in V \text{ such that } xy \in E\}.$$

Note that if we define

$$(3.2) \quad \Omega_G \doteq \{\omega \in S^V : \omega_x + \omega_y \leq C \text{ for all } x, y \in V \text{ such that } xy \in E\},$$

then, for every $\omega \in \hat{\Omega}_G$, $\hat{\sigma}_V(\omega) \in \Omega_G$.

3.1.2. *Markov specifications and Gibbs measures.* We recall here the definition of a Markov specification (see, e.g., [14, 28]), adapted, where necessary, to the context of this application.

DEFINITION 3.1. A Markov specification on $(\hat{S}^{\hat{V}}, \hat{\mathcal{F}}(\hat{\mathcal{V}}))$ is a collection $\hat{\Pi} = \{\hat{\pi}_U\}_{U \in \hat{\mathcal{V}}}$ of stochastic kernels $\hat{\pi}_U : \hat{S}^{\hat{\partial} U} \times \hat{\mathcal{F}}(U) \to [0, 1]$ that satisfy the following consistency condition: for each $U, W \in \hat{\mathcal{V}}$ with $U \subset W$ and $\tau \in \hat{S}^{\hat{V}}$,

$$\hat{\pi}_W(\tau_{\hat{\partial} W}; \hat{\sigma}_W = \tau_W) = \frac{\hat{\pi}_U(\tau_{\hat{\partial} U}; \hat{\sigma}_U = \tau_U)}{\hat{h}_{W,U}(\tau_{\Delta W \setminus U})}$$

for some normalizing function $\hat{h}_{W,U} : \hat{S}^{\Delta W \setminus U} \to \mathbb{R}_+$.

Recall that, by the definition of a stochastic kernel, for each $U \in \hat{\mathcal{V}}$ and $\tau \in \hat{S}^{\hat{V}}$, $\hat{\pi}_U(\tau_{\hat{\partial} U}; \cdot)$ is a probability measure on $(\hat{S}^{\hat{V}}, \hat{\mathcal{F}}(\hat{V}))$.

A generic way of obtaining Markov specifications on an arbitrary graph $G$ is to consider the collection of regular conditional probabilities associated with a consistent family of "product-form" probability distributions defined on all finite subgraphs of $G$. Specifically, given (weight) vectors $\vec{\nu} \in \mathbb{R}^{C_V + 1}$ and $\vec{\lambda} \in \mathbb{R}^{C_E + 1}$, for any $U \in \hat{\mathcal{V}}$, first let the associated distribution $\hat{m}_U \doteq \hat{m}_U^{\vec{\nu}, \vec{\lambda}}$ be defined on $\hat{S}^U$ by

$$(3.3) \qquad \hat{m}_U(\hat{\sigma} = \tau) = \begin{cases} \frac{\prod_{v \in U_V} \nu_{\tau_v} \prod_{e \in U_E} \lambda_{\tau_e}}{Z_U}, & \text{if } \tau \in \hat{\Omega}_{G[U]}, \\ 0, & \text{otherwise,} \end{cases}$$

where $Z_U$ is the usual normalizing constant that makes $\hat{m}_U$ a probability measure and $\hat{\Omega}_{G[U]}$ is defined by (3.1). Note that this measure is precisely the stationary distribution on the finite graph $\hat{G}[U]$ of the controlled unicast-multicast model associated with weight vectors $\vec{\lambda}$ and $\vec{\nu}$, as described in Section 1.1 and, in particular, for the weights given in (1.4). Then, for $U \in \hat{\mathcal{V}}$, let $\hat{\pi}_U : \hat{S}^{\hat{\partial} U} \times \hat{\mathcal{F}}(U)$ be the stochastic kernel defined by

$$(3.4) \qquad \hat{\pi}_U(\tau_{\hat{\partial} U}; \hat{\sigma}_W = \tau_W) \doteq \hat{m}_{\Delta U}(\hat{\sigma}_W = \tau_W | \hat{\sigma}_{\hat{\partial} U} = \tau_{\hat{\partial} U})$$

for every $\tau_{\hat{\partial} U} \in \hat{S}^{\hat{\partial} U}$ with $\hat{m}_{\Delta U}(\hat{\sigma}_{\hat{\partial} U} = \tau_{\hat{\partial} U}) > 0$, $W \subseteq U$ and $\tau_W \in \hat{S}^W$. It is straightforward to check that this defines a consistent family of stochastic kernels in the sense of Definition 3.1.

While the expression in (3.3) is well defined for finite sets $U$, it does not make sense for infinite sets, since $Z_U$ would be infinite. The natural generalization of the stationary distributions $\{\hat{m}_U\}_{U \in \hat{\mathcal{V}}}$ to infinite graphs turns out to be Gibbs measures that correspond to the Markov specification $\hat{\Pi} = \{\hat{\pi}_U\}_{U \in \mathcal{V}}$, with $\hat{\pi}_U$ given by (3.4).

DEFINITION 3.2. Given a Markov specification $\hat{\Pi}$ on $(\hat{S}^{\hat{V}}, \hat{\mathcal{F}}(\hat{V}))$, a Gibbs measure corresponding to $\hat{\Pi}$ is a probability measure $\hat{\mu}$ on $(\hat{S}^{\hat{V}}, \hat{\mathcal{F}}(\hat{V}))$ that satisfies, for all $U \in \hat{\mathcal{V}}$ and $\hat{\mu}$ a.s. $\tau \in \hat{S}^{\hat{V}}$,

$$(3.5) \qquad \hat{\mu}(\hat{\sigma}_U = \tau_U | \hat{\sigma}_{\hat{V} \setminus U} = \tau_{\hat{V} \setminus U}) = \hat{\pi}_U(\tau_{\partial U}; \hat{\sigma}_U = \tau_U).$$

A Gibbs measure $\hat{\mu}$ that corresponds to the particular Markov specification $\hat{\Pi} = \{\hat{\pi}_U\}_{U \in \hat{\mathcal{V}}}$ defined by (3.4) will be called a Gibbs measure for the controlled (unicast-multicast) model. We let $\mathcal{G}(\hat{\Pi})$ denote the space of Gibbs measures corresponding to the Markov specification $\hat{\Pi}$.

In Section 3.2 (see Corollary 3.7) we establish a one-to-one correspondence between Gibbs measures for the controlled model and certain Gibbs measures on $(S^V, \mathcal{F}(V))$. For this, we need to introduce analogous definitions of Markov specifications and Gibbs measures on $S^V$.

DEFINITION 3.3. A Markov specification on $(S^V, \mathcal{F}(V))$ is a collection $\Pi = \{\pi_U\}_{U \in \mathcal{V}}$ of stochastic kernels $\pi_U : S^{\partial U} \times \mathcal{F}(U) \to [0, 1]$ that satisfy the following consistency condition: for each $U, W \in \mathcal{V}$ with $U \subset W$, and $\tau \in S^V$:

$$\pi_W(\tau_{\partial W}; \sigma_W = \tau_W) = \frac{\pi_U(\tau_{\partial U}; \sigma_U = \tau_U)}{h_{W,U}(\tau_{\Delta W \setminus U})}$$

for some normalizing function $h_{W,U} : S^{\Delta W \setminus U} \to \mathbb{R}_+$.

DEFINITION 3.4. Given a Markov specification $\Pi = \{\pi_U\}_{U \in \mathcal{V}}$ on $(S^V, \mathcal{F}(V))$, a Gibbs measure corresponding to $\Pi$ is a probability measure $\mu$ on $(S^V, \mathcal{F}(V))$ that satisfies, for all $U \in \mathcal{V}$ and $\mu$ a.s. $\tau \in S^V$,

$$(3.6) \qquad \mu(\sigma_U = \tau_U | \sigma_{V \setminus U} = \tau_{V \setminus U}) = \pi_U(\tau_{\partial U}; \sigma_U = \tau_U).$$

Analogously to Definition 3.2, we will denote by $\mathcal{G}(\Pi)$ the space of Gibbs measures corresponding to the Markov specification $\Pi$.

It turns out that the Markov specifications on $(S^V, \mathcal{F}(V))$ that arise in our analysis cannot be generated from a family of consistent product-form probability distributions. Instead, we will require a more general way of generating a Markov specification on $(S^V, \mathcal{F}(V))$, namely, via an interaction function. Specifically, a mapping $\phi : S \times S \to \mathbb{R}_+$ is said to be an *interaction function* if and only if there exists some reference element $k_0 \in S$ such that

$$\begin{aligned} \phi(i, k_0) &> 0, & i &\in S, \\ \phi(i, j) &= \phi(j, i), & i, j &\in S. \end{aligned}$$

(We have omitted a third condition on the mapping $\phi$ given in [28] that is needed to ensure that the associated Markov specification is bounded, since it is automatically satisfied here, due to the finiteness of $S$.) Given an interaction function $\phi$, it was shown in [10, 11] that the collection $\{\pi_{\phi,U}\}_{U \in \mathcal{V}}$ given by

$$(3.7) \qquad \pi_{\phi,U}(\tau_{\partial U}; \sigma_U = \tau_U) \doteq \frac{\prod_{uv \in \Delta U \cap E} \phi(\tau_u, \tau_v)}{k_{\phi,U}(\tau_{\partial U})}$$

defines a consistent family of stochastic kernels $\pi_{\phi,U} \colon S^{\partial U} \times \mathcal{F}(U) \to \mathbb{R}_+$, in the sense of Definition 3.3. Here $k_{\phi,U} \colon S^{\partial U} \to \mathbb{R}_+$ is the appropriate normalizing mapping that, for each $\tau_{\partial U} \in S^{\partial U}$, makes $\pi_{\phi,U}(\tau_{\partial U}; \cdot)$ a probability measure. Given an interaction function $\phi$ and the collection $\pi_{\phi,U}$, $U \in \mathcal{V}$, defined by (3.7), we use $\Pi_\phi = \{\pi_{\phi,U}\}_{U \in \mathcal{V}}$ to represent the Markov specification with interaction function $\phi$.

3.2. *Gibbs measures for the controlled model.* The first main result of this section is Theorem 3.1 which, for a given pair of weight vectors $\vec{\nu}$ and $\vec{\lambda}$, establishes a one-to-one correspondence between the set $\mathcal{G}(\hat{\Pi})$ of Gibbs measures for the controlled model and the set $\mathcal{G}(\Pi_\phi)$ of Gibbs measures corresponding to the Markov specification $\Pi_\phi$ with interaction function $\phi$ defined by

$$(3.8) \qquad \phi(i,j) \doteq \begin{cases} (\nu_i \nu_j)^{1/(q+1)} \displaystyle\sum_{k=0}^{\min(C_E, C-(i+j))} \lambda_k, \\ \qquad \text{for } i,j \in S \colon 0 \le i+j \le C, \\ 0, \qquad \text{otherwise.} \end{cases}$$

Note that this function $\phi$ is clearly symmetric and 0 serves as a reference element, so that $\phi$ is indeed an interaction function. In fact, it gives rise to a repulsive Markov specification, in the sense defined by Zachary [29], since, for all $i \le k, j \le l$,

$$\frac{\phi(i,j)}{\phi(i,l)} \le \frac{\phi(k,j)}{\phi(k,l)}.$$

The proof of the theorem relies on the following two lemmas. The first lemma is an elementary graph-theoretic result, while the second lemma establishes a relation between the specifications $\hat{\Pi}$ and $\Pi_\phi$.

LEMMA 3.5. *For any $U \subseteq V$, $\hat{\partial}[U \cup [E \cap \hat{\partial}U]] = \partial U$. Moreover, for any $U \subset \hat{V}$, $\hat{\partial}(\Delta U \cap E) = \Delta U \cap V$.*

PROOF. First, note that for $U \subset V \subset \hat{V}$, by the definition of $\partial$ and $\hat{\partial}$, it follows that $\hat{\partial}U = [E \cap \hat{\partial}U] \cup \partial U$ and $\hat{\partial}[E \cap \hat{\partial}U] \subseteq \partial U \cup U$. The first relation in the lemma then follows from the observation that

$$\hat{\partial}[A \cup B] = [\hat{\partial}A \cup \hat{\partial}B] \setminus [A \cup B] \qquad \text{for all } A, B \subset \hat{V}.$$

For the proof of the second relation, first note that, by the last display and the fact that $A \subset E$ implies $\hat{\partial} A \subset E$, we have $\hat{\partial}[\Delta U \cap E] = \hat{\partial}[U_E \cup [\hat{\partial} U_V \cap E]] = \hat{\partial} U_E \cup \hat{\partial}[\hat{\partial} U_V \cap E]$ and $\Delta U \cap V = \hat{\partial} U_E \cup \overline{U}_V$. From the last two equalities, it is clear that to complete the proof, it suffices to show that $\overline{U}_V = \hat{\partial}[\hat{\partial} U_V \cap E]$. Now, suppose $x \in \overline{U}_V$. Then, there must exist $z \in \overline{U}_V$ such that $xz \in \hat{\partial} U_V \cap E$, from which it directly follows that $x \in \hat{\partial}[\hat{\partial} U_V \cap E]$. On the other hand, if $x \in \hat{\partial}[\hat{\partial} U_V \cap E]$, then there must exist $y \in \hat{\partial} U_V \cap E$ such that $y = xz$ for some $z \in V$. However, the fact that $y \in \hat{\partial} U_V \cap E$ implies that $y = uv$ for some $u \in U_V$ and $v \in \overline{U}_V$. Since the representation for $y$ is unique, this implies that $x = u$ or $x = v$, from which it follows that $x \in \overline{U}_V$, and the proof is complete. $\square$

LEMMA 3.6. *Let $\hat{\Pi}$ be the Markov specification for the controlled unicast-multicast model and let $\Pi_\phi$ be the Markov specification associated with the interaction function $\phi$ defined by (3.8). Then, for any $\tau \in S^V$ and $U \subset V$,*

$$\hat{\pi}_{U \cup [E \cap \hat{\partial} U]}(\tau_{\partial U}; \hat{\sigma}_U = \tau_U) = \pi_{\phi, U}(\tau_{\partial U}; \sigma_U = \tau_U).$$

PROOF. Fix $\tau \in S^V$ and let $U \subset V$. Then, by Lemma 3.5, $\partial U = \hat{\partial}[U \cup [E \cap \hat{\partial} U]]$, and so $\hat{\pi}_{U \cup [E \cap \hat{\partial} U]}(\tau_{\partial U}; \hat{\sigma}_U = \tau_U) = \hat{\pi}_{U \cup [E \cap \hat{\partial} U]}(\tau_{\hat{\partial}[U \cup [E \cap \hat{\partial} U]]}; \hat{\sigma}_U = \tau_U)$ which, by (3.4) and the fact that $\Delta[U \cup [E \cap \hat{\partial} U]] = U \cup [E \cap \hat{\partial} U] \cup \partial U = \Delta U$, is equal to

$$\hat{m}_{\Delta U}(\hat{\sigma}_U = \tau_U | \hat{\sigma}_{\partial U} = \tau_{\partial U})$$

$$= \sum_{\eta \in \hat{S}^{\hat{\partial} U \cap E}} \hat{m}_{\Delta U}(\hat{\sigma}_{U \cup [\hat{\partial} U \cap E]} = \tau_U \eta | \hat{\sigma}_{\partial U} = \tau_{\partial U})$$

$$= \frac{\sum_{\eta \in \hat{S}^{\hat{\partial} U \cap E}} \hat{m}_{\Delta U}(\hat{\sigma}_{\Delta U} = \tau_{\overline{U}} \eta)}{\hat{m}_{\Delta U}(\hat{\sigma}_{\partial U} = \tau_{\partial U})}$$

$$= \frac{\prod_{u \in U} \nu_{\tau_u} \prod_{v \in \partial U} \nu_{\tau_v} \prod_{xy \in \Delta U \cap E} \sum_{k=0}^{\min(C_E, C - \tau_x - \tau_y)} \lambda_k}{\hat{m}_{\Delta U}(\hat{\sigma}_{\partial U} = \tau_{\partial U})}$$

if $0 \leq \tau_u + \tau_v \leq C$ for $uv \in E \cap \Delta U$, and 0 otherwise. For $v \in \partial U$, let $r_v \doteq |\{uv : u \in U, v \in \partial U\}|$ be the number of neighbors of $v$ in $U$ (clearly $r_v \leq q+1$). Then, each $u \in U$ has $q+1$ edges in $\Delta U$ emanating from it and each $v \in \partial U$ has $r_v$ edges in $\Delta U$ emanating from it. Consequently, when $0 \leq \tau_u + \tau_v \leq C$ for every $uv \in E \cap \Delta U$, the right-hand side of the last display can be rewritten as

$$\frac{(\prod_{v \in \partial U} \nu_{\tau_v}^{(q+1-r_v)/(q+1)}) \prod_{xy \in \Delta U \cap E} (\nu_{\tau_x} \nu_{\tau_y})^{1/(q+1)} \sum_{k=0}^{\min(C_E, C - \tau_x - \tau_y)} \lambda_k}{m_{\Delta U}(\hat{\sigma}_{\partial U} = \tau_{\partial U})}.$$

Now, define

$$k_{\phi,U}(\tau_{\partial U}) \doteq \left( \prod_{v \in \partial U} \nu_{\tau_v}^{(r_v - q - 1)/(q+1)} \right) m_{\Delta U}(\hat{\sigma}_{\partial U} = \tau_{\partial U}).$$

Then, combining the last three displays and recalling the definitions of $\phi$ and $\pi_{\phi,U}$ given in (3.8) and (3.7), respectively, we obtain

$$\hat{\pi}_{U \cup [E \cap \partial U]}(\tau_{\partial U}; \hat{\sigma}_U = \tau_U) = \frac{\prod_{xy \in \Delta U \cap E} \phi(\tau_x, \tau_y)}{k_{\phi,U}(\tau_{\partial U})} = \pi_{\phi,U}(\tau_{\partial U}; \sigma_U = \tau_U),$$

which completes the proof of the lemma. □

Given a probability measure $\hat{\mu}$ on $\hat{S}^{\hat{V}}$, the projection $\mu$ of $\hat{\mu}$ onto $V$ is defined in the obvious way: for all $\tau \in S^V$ and $U \subseteq V$,

$$\mu(\sigma_U = \tau_U) = \hat{\mu}(\omega \in \hat{S}^{\hat{V}} : \hat{\sigma}_U(\omega) = \tau_U) = \hat{\mu}(\hat{\sigma}_U = \tau_U).$$

We now state and prove Theorem 3.1.

THEOREM 3.1. *Suppose positive integers $C_V, C_E, C$ and weight vectors $\vec{\nu} \in \mathbb{R}_+^{C_V+1}$, $\vec{\lambda} \in \mathbb{R}_+^{C_E+1}$ are given. Let $\hat{\Pi}$ be the Markov specification for the associated controlled unicast-multicast model and let $\Pi_\phi$ be the Markov specification with interaction function $\phi$ defined by (3.8). Then, the following three properties hold:*

1. *Given $\hat{\mu} \in \mathcal{G}(\hat{\Pi})$, its projection $\mu$ onto $S^V$ lies in $\mathcal{G}(\Pi_\phi)$;*
2. *Any Gibbs measure $\hat{\mu} \in \mathcal{G}(\hat{\Pi})$ is uniquely determined by its projection onto $S^V$;*
3. *Given any $\mu \in \mathcal{G}(\Pi_\phi)$, there exists $\hat{\mu} \in \mathcal{G}(\hat{\Pi})$ such that the projection of $\hat{\mu}$ onto $S^V$ is equal to $\mu$.*

PROOF. Let $\hat{\mu} \in \mathcal{G}(\hat{\Pi})$ and let $\mu$ be its projection onto $S^V$. Then, for any $U \in \mathcal{V}$ and $\tau \in \Omega_G$ [recall the definition of $\Omega_G$ given in (3.2)],

$$\begin{aligned}
(3.9) \quad \mu(\sigma_U = \tau_U | \sigma_{V \setminus U} = \tau_{V \setminus U}) &= \hat{\mu}(\hat{\sigma}_V = \tau_V | \hat{\sigma}_{V \setminus U} = \tau_{V \setminus U}) \\
&= \sum_{\eta \in \hat{S}^E} \hat{\mu}(\hat{\sigma} = \tau_V \eta | \hat{\sigma}_{V \setminus U} = \tau_{V \setminus U}).
\end{aligned}$$

Given any $\eta \in \hat{S}^E$, $\hat{\mu}(\hat{\sigma} = \tau_V \eta | \hat{\sigma}_{V \setminus U} = \tau_{V \setminus U})$ is equal to

$$\begin{aligned}
\hat{\mu}(\hat{\sigma}_{U \cup [E \cap \partial U]} = \tau_U \eta_{E \cap \partial U} | \hat{\sigma}_{[V \setminus U] \cup [E \setminus \partial U]} = \tau_{V \setminus U} \eta_{E \setminus \partial U}) \\
\times \hat{\mu}(\hat{\sigma}_{E \setminus \partial U} = \eta_{E \setminus \partial U} | \hat{\sigma}_{V \setminus U} = \tau_{V \setminus U}).
\end{aligned}$$

Note that $U \cup [E \cap \partial U] \in \hat{\mathcal{V}}$, since $U \in \mathcal{V}$ and the tree is a locally finite graph. Moreover, $\hat{V} \setminus [U \cup [E \cap \partial U]] = [V \setminus U] \cup [E \setminus \partial U]$. Therefore, using the relation $\hat{\partial}[U \cup [E \cap \hat{\partial}U]] = \partial U \subseteq V \setminus U$ proved in Lemma 3.5, along with the fact that $\hat{\mu}$ satisfies (3.5), we see that

$$\hat{\mu}(\hat{\sigma}_{U \cup [E \cap \hat{\partial}U]} = \tau_U \eta_{E \cap \hat{\partial}U} | \hat{\sigma}_{[V \setminus U] \cup [E \setminus \hat{\partial}U]} = \tau_{V \setminus U} \eta_{E \setminus \hat{\partial}U})$$
$$= \hat{\pi}_{U \cup [E \cap \hat{\partial}U]}(\tau_{\partial U}; \hat{\sigma}_{U \cup [E \cap \hat{\partial}U]} = \tau_U \eta_{E \cap \hat{\partial}U}).$$

The last two displays, when substituted back into (3.9), show that for every $U \in \mathcal{V}$ and $\tau \in \Omega_G$,

$$\mu(\sigma_U = \tau_U | \sigma_{V \setminus U} = \tau_{V \setminus U}) = \sum_{\eta \in \hat{S}^E} \hat{\pi}_{U \cup [E \cap \hat{\partial}U]}(\tau_{\partial U}; \hat{\sigma}_{U \cup [E \cap \hat{\partial}U]} = \tau_U \eta_{E \cap \hat{\partial}U})$$
$$\times \hat{\mu}(\hat{\sigma}_{E \setminus \hat{\partial}U} = \eta_{E \setminus \hat{\partial}U} | \hat{\sigma}_{V \setminus U} = \tau_{V \setminus U})$$
$$= \sum_{\eta' \in \hat{S}^{E \cap \hat{\partial}U}} \hat{\pi}_{U \cup [E \cap \hat{\partial}U]}(\tau_{\partial U}; \hat{\sigma}_{U \cup [E \cap \hat{\partial}U]} = \tau_U \eta')$$
$$\times \sum_{\eta'' \in \hat{S}^{E \setminus \hat{\partial}U}} \hat{\mu}(\hat{\sigma}_{E \setminus \hat{\partial}U} = \eta'' | \hat{\sigma}_{V \setminus U} = \tau_{V \setminus U})$$
$$= \sum_{\eta' \in \hat{S}^{E \cap \hat{\partial}U}} \hat{\pi}_{U \cup [E \cap \hat{\partial}U]}(\tau_{\partial U}; \hat{\sigma}_{U \cup [E \cap \hat{\partial}U]} = \tau_U \eta')$$
$$= \hat{\pi}_{U \cup [E \cap \hat{\partial}U]}(\tau_{\partial U}; \hat{\sigma}_U = \tau_U),$$

where the second equality uses the fact that, for any $U \in \mathcal{V}$, $\hat{S}^E = \hat{S}^{E \cap \hat{\partial}U} \times \hat{S}^{E \setminus \hat{\partial}U}$, and the third equality follows because of the fact that, since $\hat{\mu}$ is a probability measure concentrated on $\hat{S}^{V \cup E}$, the second summation in the previous line is equal to 1. The last display, when combined with Lemma 3.6, then shows that for every $U \in \mathcal{V}$ and $\mu$ a.s. $\tau \in S^V$,

$$\mu(\sigma_U = \tau_U | \sigma_{V \setminus U} = \tau_{V \setminus U}) = \pi_{\phi, U}(\tau_{\partial U}; \sigma_U = \tau_U),$$

which, by Definition 3.4, implies that $\mu \in \mathcal{G}(\Pi_\phi)$, thus establishing the first property of the theorem.

To prove the second property, let $\hat{\mu} \in \mathcal{G}(\hat{\Pi})$ and let $\mu$ be its projection onto $S^V$. Then, for $U \in \hat{\mathcal{V}}$ and $\tau \in \hat{S}^{\hat{V}}$, $\hat{\mu}(\hat{\sigma}_U = \tau_U)$ is equal to

$$\sum_{\eta \in \hat{S}^{\hat{\partial}U}} \hat{\mu}(\hat{\sigma}_{\Delta U} = \tau_U \eta) = \sum_{\eta \in \hat{S}^{\hat{\partial}U}} \hat{\mu}(\hat{\sigma}_{\Delta U \cap E} = \tau_{U_E} \eta_{\hat{\partial}U \cap E} | \hat{\sigma}_{\Delta U \cap V} = \tau_{U_V} \eta_{\hat{\partial}U \cap V})$$
$$\times \hat{\mu}(\hat{\sigma}_{\Delta U \cap V} = \tau_{U_V} \eta_{\hat{\partial}U \cap V})$$
$$= \sum_{\eta \in \hat{S}^{\hat{\partial}U}} \hat{\mu}(\hat{\sigma}_{\Delta U \cap E} = \tau_{U_E} \eta_{\hat{\partial}U \cap E} | \hat{\sigma}_{\Delta U \cap V} = \tau_{U_V} \eta_{\hat{\partial}U \cap V})$$
$$\times \mu(\sigma_{\Delta U \cap V} = \tau_{U_V} \eta_{\hat{\partial}U \cap V}),$$

where the last equality uses the fact that $\mu$ is a projection of $\hat{\mu}$ onto $S^V$. Combining the above display with the fact that $\hat{\mu}$ satisfies (3.5), and using the relation $\hat{\partial}(\Delta U \cap E) = \Delta U \cap V$ proved in Lemma 3.5, we see that

$$\hat{\mu}(\hat{\sigma}_U = \tau_U) = \sum_{\eta \in \hat{S}^{\partial U}} \hat{\pi}_{\Delta U \cap E}(\tau_{U_V} \eta_{\hat{\partial} U \cap V}; \hat{\sigma}_{\Delta U \cap E} = \tau_{U_E} \eta_{\hat{\partial} U \cap E})$$

$$\times \mu(\sigma_{\Delta U \cap V} = \tau_{V_V} \eta_{\hat{\partial} U \cap V}).$$

Since the last relation holds for every $U \in \hat{\mathcal{V}}$, given $\mu$, it uniquely characterizes $\hat{\mu}$.

For the last property, suppose $\mu \in \mathcal{G}(\Pi_\phi)$. Let $\theta$ be the probability distribution on $S_{C_E}$ given by

$$\theta(i) = \frac{\lambda_i}{\sum_{i=0}^{C_E} \lambda_i} \qquad \text{for } i = 0, 1, \dots, C_E.$$

Then, define $\hat{\mu}$ as follows: to generate a sample $\hat{\sigma}$ from $\hat{\mu}$, first pick a configuration $\sigma$ according to $\mu$, and then, independently on each edge $xy$ of $E$, choose a configuration $\hat{\sigma}_{xy}$ on $S_{C_E}$ according to the distribution $\theta$ conditioned on the constraint that $\hat{\sigma}_{xy} + \sigma_x + \sigma_y \leq C$. As a direct consequence of this definition, we see that $\hat{\mu}$ is supported on $\hat{\Omega}_{\hat{G}}$ and the projection of $\hat{\mu}$ onto $S^V$ is equal to $\mu$. Using arguments similar to those given above to prove properties 1 and 2, it is easy to verify that, in fact, $\hat{\mu} \in \mathcal{G}(\hat{\Pi})$. We leave the details to the reader.  $\square$

COROLLARY 3.7.  *Let $\hat{\Pi}$ and $\Pi_\phi$ be as defined in Theorem 3.1. Then, there is a one-to-one correspondence between the sets $\mathcal{G}(\hat{\Pi})$ and $\mathcal{G}(\Pi_\phi)$.*

PROOF.  Let $\chi$ be the mapping that takes $\mu \in \mathcal{G}(\hat{\Pi})$ to its projection onto $S^V$. Then, property 1 of Theorem 3.1 shows that $\chi$ maps $\mathcal{G}(\hat{\Pi})$ to $\mathcal{G}(\Pi_\phi)$, property 2 shows that $\chi$ is one-to-one and property 3 shows that the mapping is onto. Thus, $\chi$ defines a one-to-one correspondence between the two sets of Gibbs measures.  $\square$

As a consequence of Corollary 3.7, the analysis of phase transitions for Gibbs measures that correspond to the simple product-form specification $\hat{\Pi}$ that is Markov with respect to the more complicated graph $\hat{G}$ is reduced to the analysis of Gibbs measures corresponding to the more complicated specification $\Pi_\phi$, which is Markov with respect to the simpler tree graph $G$. The advantage of this reduction is that Gibbs measures that are Markov with respect to trees are usually much easier to analyze than Gibbs measures that are Markov with respect to more complicated graphs. However, it is worth emphasizing that, while the original Markov specification $\hat{\Pi}$ is of product-form, the associated Markov specification $\Pi_\phi$ is no longer of product-form. In addition, the set of allowable configurations in $S^V$ cannot be represented as the set of homomorphisms from the tree to a constraint graph, as was done in [5]. Hence, this problem does not fall directly within the framework considered in [5]. Instead, we use the full generality of the results in [28, 29], which consider Gibbs measures associated with general Markov specifications on trees, to prove Theorem 1.1, which shows that the analysis of Gibbs measures for the controlled unicast-multicast model can be further reduced to the study of recursions of the random field map $\Phi^{\nu,\lambda}$ given in (1.6). These recursions are analyzed for the case $C_V = 1$ in Section 4.

PROOF OF THEOREM 1.1. Let $\Psi_\phi$ be the set of vectors $\boldsymbol{\psi} \in \mathbb{R}_+^{C_V+1}$ having first component $\boldsymbol{\psi}_0 = 1$. Let the mapping $F^\phi : \Psi_\phi \to \Psi_\phi$ be given by

$$F_i^\phi(\boldsymbol{\psi}) = \left( \frac{\sum_{j=0}^{\min(C-i,C_V)} \phi(i,j)\psi_j}{\sum_{j=0}^{\min(C,C_V)} \phi(0,j)\psi_j} \right)^q \qquad \text{for } i = 0, 1, \dots, C_V.$$

Define $\boldsymbol{\psi}^{(0)} \in \mathbb{R}_+^{C_V+1}$ to be the vector whose first component is equal to 1 and remaining components are equal to zero, that is, $\boldsymbol{\psi}^{(0)} = (1, 0, \dots, 0)$, and let $\boldsymbol{\psi}^{(n+1)} = F^\phi(\boldsymbol{\psi}^{(n)})$ be the iterates of $\boldsymbol{\psi}^{(0)}$ under the map $F^\phi$. Then, Theorem 4.1 of [29] states that the Markov specification $\Pi_\phi$ has a unique Gibbs measure if and only if the sequence of iterates $\boldsymbol{\psi}^{(n)}, n \geq 0$, converges to a limit and, in this case, the limit $\boldsymbol{\psi}^*$ must be the unique fixed point of the map $F^\phi$.

Consider the map $A : \Psi_\phi \to \mathbb{R}_+^{C_V}$ defined by

$$A_k(\boldsymbol{\psi}) = \nu_k^{1/(q+1)} \psi_k \qquad \text{for } k = 1, \dots, C_V,$$

and let $A^{-1} : \mathbb{R}_+^{C_V} \to \Psi_\phi$ be the map given by

$$A_k^{-1}(\xi) = \begin{cases} 1, & \text{for } k = 0, \\ \nu_k^{-1/(q+1)} \xi_k, & \text{for } k = 1, \dots, C_V, \end{cases}$$

where $A_k(\phi)$ is used to denote the $k$th component of $A(\phi)$. Then, clearly, $A \circ A^{-1}$ and $A^{-1} \circ A$ are the identity mappings on $\mathbb{R}_+^{C_V}$ and $\Psi_\phi$, respectively. We now claim that, for any $\boldsymbol{\psi} \in \Psi_\phi$,

$$(3.10) \qquad F^\phi(\boldsymbol{\psi}) = A^{-1} \Phi^{\vec{\nu}, \vec{\lambda}}(A\boldsymbol{\psi}).$$

It is easy to see that the theorem follows from this claim. Indeed, let $\mathbf{0}$ be the zero vector in $\mathbb{R}_+^{C_V}$ and note that $\mathbf{0} = A(\boldsymbol{\psi}^{(0)})$. Moreover, the relation (3.10) shows that $\boldsymbol{\psi}^*$ is a fixed point of $F^\phi$ if and only if $A(\boldsymbol{\psi}^*)$ is a fixed point of $\Phi^{\vec{\nu}, \vec{\lambda}}$ and, by induction, that $A[F^\phi]^n(\boldsymbol{\psi}) = [\Phi^{\vec{\nu}, \vec{\lambda}}]^n(A\boldsymbol{\psi})$. From this, it immediately follows that $[F^\phi]^n(\psi^{(0)}) \to \psi^*$ if and only if $[\Phi^{\vec{\nu},\vec{\lambda}}]^n(\mathbf{0}) \to A(\psi^*)$. With this equivalence, the theorem is a direct consequence of the result from Theorem 4.1 of [29] quoted above.

So it only remains to prove the claim (3.10). Suppose $\psi \in \Psi_\phi$. Then, substituting the definition of $\phi$ from (3.8), we have, for $i = 1, \ldots, C_V$,

$$
\begin{aligned}
F_i^\phi(\psi) &= \left( \frac{\sum_{j=0}^{\min(C-i,C_V)} \phi(i,j)\psi_j}{\sum_{j=0}^{\min(C,C_V)} \phi(0,j)\psi_j} \right)^q \\
&= \left( \frac{\sum_{j=0}^{\min(C-i,C_V)} (\nu_i \nu_j)^{1/(q+1)} \sum_{k=0}^{\min(C-(i+j),C_E)} \lambda_k \psi_j}{\sum_{j=0}^{\min(C,C_V)} (\nu_0 \nu_j)^{1/(q+1)} \sum_{k=0}^{\min(C-j,C_E)} \lambda_k \psi_j} \right)^q \\
&= \left( \frac{\sum_{k=0}^{\min(C-i,C_E)} \lambda_k \sum_{j=0}^{\min(C-(i+k),C_V)} (\nu_i \nu_j)^{1/(q+1)} \psi_j}{\sum_{k=0}^{\min(C,C_E)} \lambda_k \sum_{j=0}^{\min(C-k,C_V)} \nu_j^{1/(q+1)} \psi_j} \right)^q \\
&= \nu_i^{q/(q+1)} \left( \frac{\sum_{k=0}^{\min(C-i,C_E)} \lambda_k \sum_{j=0}^{\min(C-(i+k),C_V)} \nu_j^{1/(q+1)} \psi_j}{\sum_{k=0}^{C_E} \lambda_k \sum_{j=0}^{\min(C-k,C_V)} \nu_j^{1/(q+1)} \psi_j} \right)^q \\
&= \nu_i^{-1/(q+1)} \nu_i \left( \frac{\sum_{k=0}^{\min(C-i,C_E)} \lambda_k [1 + \sum_{j=1}^{\min(C-(i+k),C_V)} A_j(\psi)]}{\sum_{k=0}^{C_E} \lambda_k [1 + \sum_{j=1}^{\min(C-k,C_V)} A_j(\psi)]} \right)^q \\
&= \nu_i^{-1/(q+1)} \Phi_i^{\vec{\nu},\vec{\lambda}}(A(\psi)),
\end{aligned}
$$

where the third equality is obtained by exchanging the order of summation, the fourth equality uses the fact that $C_E \leq C$ (which holds by assumption) to replace $\min(C_E, C)$ by $C_E$ and the fifth equality holds because, for every $k = 0, \ldots, C-i$, $\min(C-i-k, C_V) \geq 0$ and, likewise, for every $k = 0, \ldots, C$, $\min(C-k, C_V) \geq 0$. This establishes (3.10) and, therefore, completes the proof of the theorem. $\quad\square$

## 4. Analysis for the case $C_V = 1$, $C_E = C$.

In Section 3 we established a relation between Gibbs measures for the controlled unicast-multicast model on the $(q+1)$-regular tree and recursions of the random field map $\Phi^{\vec{\nu},\vec{\lambda}}$. In this section we analyze this recursion for the case $C_V = 1$ and $C_E$ equal to some integer $C \geq 2$. (Note that the restriction $C \geq 2$ is without loss of generality, since, as mentioned earlier, when $C_E = C_V = 1$, the model is equivalent to a pure multicast model with a modified weight, which is well understood [18, 27, 28, 29].) When $C_V = 1$, $\vec{\nu} = (1, \nu)$ can be parametrized by just the scalar $\nu$, and $\Phi^{\vec{\nu},\vec{\lambda}} = \Phi^{\nu,\vec{\lambda}} : \mathbb{R}_+ \to \mathbb{R}_+$ reduces to the following

one-dimensional map:

$$\Phi^{\nu,\vec{\lambda}}(\xi) = \nu\left(\frac{\sum_{i=0}^{C-1} \lambda_i[1 + \sum_{j=1}^{\min(C-1-i,1)} \xi]}{\sum_{i=0}^{C} \lambda_i[1 + \sum_{j=1}^{\min(C-i,1)} \xi]}\right)^q$$

$$= \nu\left(\frac{\sum_{i=0}^{C-1} \lambda_i + \xi\sum_{i=0}^{C-2} \lambda_i}{\sum_{i=0}^{C} \lambda_i + \xi\sum_{i=0}^{C-1} \lambda_i}\right)^q$$

$$= \nu\left(\frac{\Lambda_{C-1} + \xi\Lambda_{C-2}}{\Lambda_C + \xi\Lambda_{C-1}}\right)^q,$$

where one should recall that $\Lambda_n \doteq \sum_{i=0}^{n} \lambda_i$. Throughout this section, we will assume that the weight vector $\vec{\lambda}$ satisfies condition (1.7) stated in Assumption 1.1. In Lemma 4.1, we now verify this condition for the weight vectors defined in (1.4) and (1.5).

LEMMA 4.1. *Given any positive scalar $\lambda$, Assumption 1.1 is satisfied by the weight vectors $\vec{\lambda}$ defined in terms of $\lambda$ by (1.4) and (1.5).*

PROOF. Fix $C \geq 2$. We first consider the case of the "Poisson" weight vector defined by (1.4), with $\lambda_i = \lambda^i/i!$. In this case,

$$\Lambda_{C-1}^2 - \Lambda_C\Lambda_{C-2} = \Lambda_{C-1}\left(\Lambda_{C-2} + \frac{\lambda^{C-1}}{(C-1)!}\right) - \left(\Lambda_{C-1} + \frac{\lambda^C}{C!}\right)\Lambda_{C-2}$$

$$= \frac{\lambda^{C-1}}{(C-1)!}\Lambda_{C-1} - \frac{\lambda^C}{C!}\Lambda_{C-2}$$

$$= \frac{\lambda^{C-1}}{(C-1)!} + \sum_{i=0}^{C-2}\left(\frac{\lambda^{C-1}}{(C-1)!}\frac{\lambda^{i+1}}{(i+1)!} - \frac{\lambda^C}{C!}\frac{\lambda^i}{i!}\right).$$

Now, the coefficient for $\lambda^{C+i}$ is positive for $0 \leq i \leq C-2$ if and only if

$$\frac{1}{(C-1)!}\frac{1}{(i+1)!} - \frac{1}{C!}\frac{1}{i!} > 0,$$

and this reduces to the condition that $C > i+1$, which certainly holds for $0 \leq i \leq C-2$. Thus, each term is greater than zero, and (1.7) follows.

Now, suppose the weight vector is defined by (1.5), so that $\lambda_i = \lambda^i$. Then, simple algebraic manipulations show that

$$(4.1) \qquad \Lambda_{C-1}^2 - \Lambda_C\Lambda_{C-2} = \lambda^{C-1},$$

which is positive, since $\lambda > 0$. Thus, once again, (1.7) is established. □

In Section 4.1 we fully characterize the phase transition surface (namely the boundary in parameter space between the regions of uniqueness and multiplicity of Gibbs measures) for this model and show that the phase transition, whenever it occurs, is nonmonotone in the parameter $\nu$ (for fixed $q, C$ and $\vec{\lambda}$). In Section 4.2 we study monotonicity of the phase transitions with respect to the other parameters $q$ and $\lambda$.

4.1. *Characterization of the phase transition surface.* In Lemma 4.2 we show that the map $\Phi^{\nu,\vec{\lambda}}$ always has a unique fixed point.

LEMMA 4.2. *Suppose $\vec{\lambda}$ satisfies Assumption 1.1 and $\nu > 0$. Then the following two properties hold:*

1. *The map $\Phi^{\nu,\vec{\lambda}}$ has exactly one fixed point, $\xi^* = \xi^*(\nu, \vec{\lambda})$, in the interval $[0, \infty)$.*
2. *For a fixed $\vec{\lambda}$, the range $\{\xi^*(\nu, \vec{\lambda}), \nu \in (0, \infty)\}$ is equal to $(0, \infty)$. Moreover, there is a one-to-one correspondence between $\nu$ and the corresponding fixed point $\xi^*$, given by $\nu = J(\xi^*)$, where $J$ is the mapping defined in (1.8).*

PROOF. Observe that $\partial \Phi^{\nu,\vec{\lambda}}(\xi)/\partial \xi$ is equal to

$$\nu q \left( \frac{\Lambda_{C-1} + \xi \Lambda_{C-2}}{\Lambda_C + \xi \Lambda_{C-1}} \right)^{q-1} \frac{(\Lambda_C + \xi \Lambda_{C-1})\Lambda_{C-2} - (\Lambda_{C-1} + \xi \Lambda_{C-2})\Lambda_{C-1}}{(\Lambda_C + \xi \Lambda_{C-1})^2}$$

$$= \nu q \left( \frac{\Lambda_{C-1} + \xi \Lambda_{C-2}}{\Lambda_C + \xi \Lambda_{C-1}} \right)^q \frac{\Lambda_C \Lambda_{C-2} - \Lambda_{C-1}^2}{(\Lambda_{C-1} + \xi \Lambda_{C-2})(\Lambda_C + \xi \Lambda_{C-1})},$$

which, in turn, shows that

$$(4.2) \qquad \frac{\partial \Phi^{\nu,\vec{\lambda}}(\xi)}{\partial \xi} = \Phi^{\nu,\vec{\lambda}}(\xi) \frac{q(\Lambda_C \Lambda_{C-2} - \Lambda_{C-1}^2)}{(\Lambda_{C-1} + \xi \Lambda_{C-2})(\Lambda_C + \xi \Lambda_{C-1})}.$$

Due to (1.7) and the fact that $\Lambda_k$, $\Phi^{\nu,\vec{\lambda}}(\xi) > 0$, the last expression is negative. Furthermore, $\Phi^{\nu,\vec{\lambda}}(0) = \nu(\Lambda_{C-1}/\Lambda_C)^q > 0$ and $\lim_{\xi \to \infty} \Phi^{\nu,\vec{\lambda}}(\xi) = \nu(\Lambda_{C-2}/\Lambda_{C-1})^q$. Hence, since $\Phi^{\nu,\vec{\lambda}}(\xi)$ is continuous and strictly decreasing in $\xi$, with $\Phi^{\nu,\vec{\lambda}}(\xi) > 0$ for all $\xi \geq 0$, $\Phi^{\nu,\vec{\lambda}}$ has a unique fixed point.

To show the second property, fix $\vec{\lambda}$ and note that $\nu$ satisfies the implicit equation

$$\nu = J(\xi^*(\nu)),$$

where $\xi^* = \xi^*(\nu)$ is the unique fixed point of $\Phi^{\nu,\vec{\lambda}}$. Evaluating the derivative of $J$ at the fixed point, we obtain

$$\frac{\partial J}{\partial \xi}(\xi^*) = \left( \frac{\Lambda_C + \xi^* \Lambda_{C-1}}{\Lambda_{C-1} + \xi^* \Lambda_{C-2}} \right)^q$$

$$+ \xi^* q \left( \frac{\Lambda_C + \xi^* \Lambda_{C-1}}{\Lambda_{C-1} + \xi^* \Lambda_{C-2}} \right)^{q-1} \frac{\Lambda_{C-1}^2 - \Lambda_C \Lambda_{C-2}}{(\Lambda_{C-1} + \xi^* \Lambda_{C-2})^2}$$

$$= \frac{\nu}{\xi^*} + q\nu \frac{\Lambda_{C-1}^2 - \Lambda_C \Lambda_{C-2}}{(\Lambda_C + \xi^* \Lambda_{C-1})(\Lambda_{C-1} + \xi^* \Lambda_{C-2})}.$$

By Assumption 1.1, $\Lambda_{C-1}^2 - \Lambda_C \Lambda_{C-2} > 0$, and so this derivative is positive. Thus, by the implicit function theorem, there is a one-to-one correspondence between $\nu$ and $\xi^*(\nu)$. Furthermore, it is easy to see from the map that $\xi^*(\nu) \to 0$ as $\nu \to 0$ and $\xi^*(\nu) \to \infty$ as $\nu \to \infty$, so that the range of $\xi^*(\nu), \nu \in (0, \infty)$ is $(0, \infty)$.   $\square$

In order to investigate when iterates of the map $\Phi^{\nu, \vec{\lambda}}$ starting at $\mathbf{0}$ converge to the unique fixed point $\xi^*$, we will need the following dynamical systems result. Given an interval $I \subset \mathbb{R}$, a fixed point $x^*$ of a map $F : I \to I$ is said to be *globally stable* if, for all $x \in I$, $F^n(x) \to x^*$ (where $F^n$ represents the $n$th iteration of the map $F$).

LEMMA 4.3. *Let $F : I \to I$ be a thrice continuously differentiable map on a bounded interval $I \subset \mathbb{R}$ that has a unique fixed point $x^*$. If $F$ has a negative Schwarzian derivative, that is, if $SF(x) < 0$ for all $x \in I$, where*

$$SF(x) \doteq \frac{d^3F/dx^3}{dF/dx} - \frac{3}{2} \left( \frac{d^2F/dx^2}{dF/dx} \right)^2,$$

*then $x^*$ is globally stable if and only if*

$$\left| \frac{dF}{dx}(x^*) \right| \leq 1.$$

*Moreover, if $F$ is decreasing on $I$, then global stability of the fixed point $x^*$ is equivalent to the condition that $F^n(l) \to x^*$, where $l$ is the left endpoint of the interval $I$.*

PROOF.   The proof of the first statement is a consequence of standard results concerning one-dimensional dynamics (see, e.g., page 158 of [9]). The fact that $F$ is monotone decreasing implies that $F \circ F$ is monotone increasing and that $F(x) \leq F(l)$ for every $x \in I$, so that we can assume, without loss of generality, that $F(l) = r$ is the right endpoint of the interval $I$. Thus, for any $x \in I$, we have

$$F^{2n}(l) \leq F^{2n}(x) \leq F^{2n}(F(l)) = F^{2n+1}(l)$$

and

$$F^{2n+1}(l) \geq F^{2n+1}(x) \geq F^{2n+2}(l),$$

which, together, show that $F^n(x) \to x^*$, for all $x \in I$ if and only if $F^n(l) \to x^*$, as desired.  □

We now show that the map $\Phi^{\nu,\vec{\lambda}}$ satisfies the conditions of Lemma 4.3.

LEMMA 4.4. *Suppose Assumption* 1.1 *is satisfied and* $\xi^* = \xi^*(\nu, \vec{\lambda})$ *is the unique fixed point of the map* $\Phi^{\nu,\vec{\lambda}}$. *Then,* $[\Phi^{\nu,\vec{\lambda}}]^n(0) \to \xi^*$ *if and only if* $|\partial \Phi^{\nu,\vec{\lambda}}(\xi^*)/\partial \xi| \leq 1$.

PROOF.  First, note that, by definition of the random field map, it follows that, for any $\xi \in [0, \infty)$, $\Phi^{\nu,\vec{\lambda}}(\xi) \in [0, \nu]$ and, therefore, we can assume, without loss of generality, that the map $\Phi^{\nu,\vec{\lambda}}$ is defined on a bounded interval $I$. As a consequence of Lemma 4.3, it is sufficient to show that $\Phi^{\nu,\vec{\lambda}}$ is thrice continuously differentiable, decreasing and has a negative Schwarzian derivative (with respect to $\xi$). It is easy to see that $\Phi^{\nu,\vec{\lambda}}$ is thrice differentiable and in fact satisfies

$$\frac{\partial \Phi^{\nu,\vec{\lambda}}(\xi)}{\partial \xi} = \Phi^{\nu,\vec{\lambda}}(\xi) \frac{q(\Lambda_C \Lambda_{C-2} - \Lambda_{C-1}^2)}{(\Lambda_{C-1} + \xi \Lambda_{C-2})(\Lambda_C + \xi \Lambda_{C-1})},$$

$$\frac{\partial^2 \Phi^{\nu,\vec{\lambda}}(\xi)}{\partial \xi^2} = \frac{\partial \Phi^{\nu,\vec{\lambda}}(\xi)}{\partial \xi} \frac{(q-1)\Lambda_C \Lambda_{C-2} - (q+1)\Lambda_{C-1}^2 - 2\xi \Lambda_{C-1}\Lambda_{C-2}}{(\Lambda_{C-1} + \xi \Lambda_{C-2})(\Lambda_C + \xi \Lambda_{C-1})},$$

$$\frac{\partial^3 \Phi^{\nu,\vec{\lambda}}(\xi)}{\partial \xi^3} = \frac{\partial \Phi^{\nu,\vec{\lambda}}(\xi)}{\partial \xi} \frac{G(\xi)}{(\Lambda_{C-1} + \xi \Lambda_{C-2})^2(\Lambda_C + \xi \Lambda_{C-1})^2},$$

where $G(\xi) = [(q-1)\Lambda_C \Lambda_{C-2} - (q+1)\Lambda_{C-1}^2 - 2\xi \Lambda_{C-1}\Lambda_{C-2}][(q-2)\Lambda_C \times \Lambda_{C-2} - (q+2)\Lambda_{C-1}^2 - 4\xi \Lambda_{C-1}\Lambda_{C-2}] - 2\Lambda_{C-1}\Lambda_{C-2}(\Lambda_{C-1} + \xi \Lambda_{C-2})(\Lambda_C + \xi \Lambda_{C-1})$. The fact that $\Phi^{\nu,\vec{\lambda}}$ is monotone decreasing then follows from the first equality above, along with Assumption 1.1. Tedious algebraic manipulations of the quantities above, also show that the Schwarzian derivative being negative reduces to the condition that

$$0 < \tfrac{1}{2}(q^2 - 1)(\Lambda_C \Lambda_{C-2} - \Lambda_{C-1}^2)^2,$$

which holds trivially.  □

We now derive some elementary consequences of Condition A introduced in Definition 1.1 that will be useful for the proof of Theorem 4.1.

LEMMA 4.5. *Suppose Assumption* 1.1 *holds. Then, the following statements are true:*

1. *If Condition* A *holds for some* $(q, C, \vec{\lambda})$, *then*

$$(4.3) \qquad (1+q)\Lambda_C\Lambda_{C-2} + (1-q)\Lambda_{C-1}^2 < 0.$$

2. *Condition* A *holds for some* $(q, C, \vec{\lambda})$ *if and only if the quadratic* $Q$ *defined in* (1.9) *has two distinct real roots* $0 < \alpha_- < \alpha_+ < \infty$, *with* $Q(\alpha) < 0$ *for* $\alpha \in (\alpha_-, \alpha_+)$.

3. *The relation in Condition* A *holds with equality if and only if* $Q$ *has a positive real double root* $\alpha_- = \alpha_+$.

Proof. For the first property, note that if Condition A holds, then

$$q(q+1)\Lambda_C\Lambda_{C-2} < (q+1)^2\Lambda_C\Lambda_{C-2} < (q-1)^2\Lambda_{C-1}^2 < q(q-1)\Lambda_{C-1}^2,$$

so that $(q+1)\Lambda_C\Lambda_{C-2} < (q-1)\Lambda_{C-1}^2$, namely (4.3), holds.

For the second property, recall that the quadratic $Q$ has two distinct real roots if and only if its discriminant is strictly positive:

$$[(1-q)\Lambda_{C-1}^2 + (1+q)\Lambda_C\Lambda_{C-2}]^2 - 4\Lambda_C\Lambda_{C-1}^2\Lambda_{C-2} > 0.$$

Thus, to establish the second property, it suffices to show that the last display is equivalent to Condition A. Indeed, using the elementary equality $(a+b)^2 - (a-b)^2 = 4ab$, the last display can be seen to be equivalent to

$$[\Lambda_C\Lambda_{C-2} + \Lambda_{C-1}^2 + q(\Lambda_C\Lambda_{C-2} - \Lambda_{C-1}^2)]^2$$
$$> (\Lambda_C\Lambda_{C-2} + \Lambda_{C-1}^2)^2 - (\Lambda_C\Lambda_{C-2} - \Lambda_{C-1}^2)^2,$$

which reduces, after straightforward algebraic manipulation, to the condition that

$$(q^2+1)(\Lambda_C\Lambda_{C-2} - \Lambda_{C-1}^2)^2 + 2q(\Lambda_C\Lambda_{C-2} + \Lambda_{C-1}^2)(\Lambda_C\Lambda_{C-2} - \Lambda_{C-1}^2) > 0.$$

By Assumption 1.1, $\Lambda_C\Lambda_{C-2} - \Lambda_{C-1}^2 < 0$, so that this reduces further to the condition that

$$(q^2+1)(\Lambda_C\Lambda_{C-2} - \Lambda_{C-1}^2) + 2q(\Lambda_C\Lambda_{C-2} + \Lambda_{C-1}^2) < 0,$$

which can be rewritten as Condition A:

$$(q+1)^2\Lambda_C\Lambda_{C-2} < (q-1)^2\Lambda_{C-1}^2.$$

The minimum of $Q$ is attained at the point

$$-[(1+q)\Lambda_C\Lambda_{C-2} + (1-q)\Lambda_{C-1}^2]/[2\Lambda_{C-1}\Lambda_{C-2}],$$

which is positive, by (4.3). Since $Q(0) > 0$, any real roots that exist must be strictly positive, with $Q(\alpha) < 0$ for $\alpha \in (\alpha_-, \alpha_+)$.

For the third property, the arguments above also show that, if the relation in Condition A holds with equality, then the discriminant is zero and the double root of the quadratic $Q$ occurs at $(q-1)\Lambda_{C-1}/\Lambda_{C-2} > 0$.  □

We are now in a position to give a precise characterization of the phase transition region. Recall the set $\mathcal{A}$ and Condition A introduced in Definition 1.1 and the function $J$ defined in (1.8). Moreover, as in Lemma 4.5, let $\alpha_-$ and $\alpha_+$ be the distinct positive roots of the quadratic $Q$ and let $\nu_- \doteq J(\alpha_-)$ and $\nu_+ \doteq J(\alpha_+)$.

THEOREM 4.1.   *Suppose Assumption 1.1 holds. Let $\xi^* = \xi^*(\nu, \vec{\lambda})$ be the unique fixed point of the map $\Phi^{\nu, \vec{\lambda}}$. Then, $|\partial \Phi(\xi^*)/\partial \xi| \leq 1$ if and only if either $(q, C, \vec{\lambda}) \notin \mathcal{A}$, or $(q, C, \vec{\lambda}) \in \mathcal{A}$ and $\nu \in (0, \nu_-] \cup [\nu_+, \infty)$.*

PROOF.   Substituting the fact that the fixed point satisfies $\xi^* = \Phi^{\nu, \vec{\lambda}}(\xi^*)$ into (4.2), we see that

$$\left| \frac{\partial \Phi^{\nu, \vec{\lambda}}(\xi^*)}{\partial \xi} \right| = \left| \frac{\xi^* q (\Lambda_C \Lambda_{C-2} - \Lambda_{C-1}^2)}{(\Lambda_{C-1} + \xi^* \Lambda_{C-2})(\Lambda_C + \xi^* \Lambda_{C-1})} \right|$$

$$= \frac{\xi^* q (\Lambda_{C-1}^2 - \Lambda_C \Lambda_{C-2})}{(\Lambda_{C-1} + \xi^* \Lambda_{C-2})(\Lambda_C + \xi^* \Lambda_{C-1})},$$

where (1.7) was used to obtain the last equality. Therefore, $|\partial \Phi(\xi^*)/\partial \xi| \leq 1$ if and only if

$$\frac{\xi^* q (\Lambda_{C-1}^2 - \Lambda_C \Lambda_{C-2})}{(\Lambda_{C-1} + \xi^* \Lambda_{C-2})(\Lambda_C + \xi^* \Lambda_{C-1})} \leq 1,$$

from which we conclude that

(4.4)                    $|\partial \Phi(\xi^*)/\partial \xi| \leq 1 \quad \Longleftrightarrow \quad Q(\xi^*) \geq 0,$

where $Q$ is the quadratic given by (1.9).

First, suppose that $(q, C, \vec{\lambda}) \notin \mathcal{A}$. Then, by Lemma 4.5, it follows that the quadratic function $Q$ has no distinct real roots. Since $Q(0) = \Lambda_C \Lambda_{C-1} > 0$, this implies that $Q(\alpha) \geq 0$, for all $\alpha \in \mathbb{R}_+$ and so, in particular, $Q(\xi^*(\nu, \vec{\lambda})) \geq 0$, for all $\nu \in (0, \infty)$. Now suppose that $(q, C, \vec{\lambda}) \in \mathcal{A}$, so that Condition A is satisfied. Then, by the second property of Lemma 4.5, it follows that $Q$ has two distinct positive real roots $\alpha_-$ and $\alpha_+$, and $Q(\alpha)$ is negative precisely when $\alpha \in (\alpha_-, \alpha_+)$. By the second property of Lemma 4.2, every $\alpha \in (0, \alpha_-] \cup [\alpha_+, \infty)$ corresponds to a fixed point of the map $\Phi^{\nu, \vec{\lambda}}$ with $\nu = J(\alpha) \in (0, J(\alpha_-)] \cup [J(\alpha_+), \infty)$ and, thus, $Q(\xi^*(\nu, \vec{\lambda})) \geq 0$ for precisely these values of $\nu$. The theorem then follows from the above arguments and (4.4).   □

The proof of Theorem 1.2 now follows.

PROOF OF THEOREM 1.2.   The proof follows as a direct consequence of Lemma 4.4, Theorems 1.1 and 4.1.   □

4.2. *Monotonicity of phase transitions with respect to* $q, \nu, \lambda$. Theorem 1.2 provides an explicit characterization of the region in the $(q, C, \vec{\lambda}, \nu)$ parameter space where there are multiple Gibbs measures. Note that, in particular, this result shows that if there is a phase transition at all for some $(q, C, \vec{\lambda}, \nu)$, then this phase transition is nonmonotone in $\nu$, in the sense that, for the same $(q, C, \vec{\lambda})$ and every sufficiently small and sufficiently large $\nu$, there exists a unique Gibbs measure. The numerical examples in Section 5 show that phase transitions do indeed occur—in the sense that for a large class of weights, the set $\mathcal{A}$ is a nontrivial subset of the parameter space—and also illustrate the nonmonotonicity with respect to $\nu$. Some intuition behind why this nonmonotonicity occurs is provided in Section 6.

We now establish some additional monotonicity properties of the phase transition. The first result, Theorem 4.2, shows that the phase transition surface is, in a sense, monotone in $q$, the degree of the tree.

THEOREM 4.2. *Suppose Assumption* 1.1 *is satisfied. If* $(q, C, \vec{\lambda}) \in \mathcal{A}$, *then* $(q + 1, C, \vec{\lambda}) \in \mathcal{A}$.

PROOF. First, note that Condition A holds for $q + 1$ if and only if $(q + 2)^2 \Lambda_C \Lambda_{C-2} < q^2 \Lambda_{C-1}^2$, which, in turn, holds if and only if

$$(4.5) \quad \begin{aligned} (q+1)^2 \Lambda_C \Lambda_{C-2} + (2q + 3)\Lambda_C \Lambda_{C-2} \\ < (q-1)^2 \Lambda_{C-1}^2 + (2q - 1)\Lambda_{C-1}^2. \end{aligned}$$

Now, suppose Condition A holds for $q$. This implies that $(q+1)^2 \Lambda_C \Lambda_{C-2} < (q-1)^2 \Lambda_{C-1}^2$ and, due to the first property of Lemma 4.5, that $(2q+2)\Lambda_C \times \Lambda_{C-2} < (2q-2)\Lambda_{C-1}^2$. In addition, due to Assumption (1.1), we also have $\Lambda_C \Lambda_{C-2} < \Lambda_{C-1}^2$. When combined, this shows that (4.5) holds and, therefore, that Condition A also holds for $q + 1$, or, equivalently, that $(q+1, C, \vec{\lambda}) \in \mathcal{A}$. □

We now show that, for the weights in (1.4), the set $\mathcal{A}$ is also monotone with respect to the parameter $\lambda$.

THEOREM 4.3. *Suppose that* $\vec{\lambda}$ *is a weight vector defined as in* (1.4) *for some* $\lambda > 0$. *Then, given any weight vector* $\vec{\lambda}'$ *defined as in* (1.4) *but with* $\lambda$ *replaced by* $\lambda' > \lambda$, $(q, C, \vec{\lambda}) \in \mathcal{A} \Rightarrow (q, C, \vec{\lambda}') \in \mathcal{A}$.

PROOF. We first show that $(q, C, \vec{\lambda}) \in \mathcal{A}$ if and only if $F(q, C, \vec{\lambda}) > 0$, where

$$F(q, C, \vec{\lambda}) \doteq (1 + q)^2 [\lambda_{C-1}\Lambda_{C-1} - \lambda_C \Lambda_{C-2}] - 4q\Lambda_{C-1}^2.$$

Indeed, by Definition 1.1, $(q, C, \vec{\lambda}) \in \mathcal{A}$ if and only if

$$(1+q)^2 \Lambda_C \Lambda_{C-2} < (q-1)^2 \Lambda_{C-1}^2$$

$$\iff (q+1)^2 (\Lambda_{C-1} + \lambda_C)(\Lambda_{C-1} - \lambda_{C-1}) < (q-1)^2 \Lambda_{C-1}^2$$

$$\iff (1+q)^2 [\lambda_{C-1} \Lambda_{C-1} - \lambda_C \Lambda_{C-2}] - 4q \Lambda_{C-1}^2 > 0.$$

Therefore, in order to prove the theorem, it suffices to show that if $F(q, C, \vec{\lambda}) > 0$, then $\partial F(q, C, \vec{\lambda})/\partial \lambda > 0$, since the latter implies that, for any $\lambda' > \lambda$, $F(q, C, \vec{\lambda}') > 0$, which, in turn, means that $(q, C, \vec{\lambda}') \in \mathcal{A}$.

Let $H(i)$ be the coefficients in the expansion of $F(q, C, \vec{\lambda})$ as a polynomial in $\lambda$: $F(q, C, \vec{\lambda}) = \sum_{i=0}^{2C-2} H(i) \lambda^i$. We now claim that the following assertion is true:

CLAIM.   *There exists $C - 1 < i^* \leq 2C - 2$ such that $H(i) \leq 0$ for all $i < i_*$ and $H(i) > 0$ for all $i \geq i^*$.*

We defer the proof of the claim and, instead, first show that the theorem follows from this claim. Indeed, if the claim holds, then for $i \geq i^*$, $iH(i)\lambda^i/i^* \geq H(i)\lambda^i > 0$ and, for $i < i^*$, $0 \geq iH(i)\lambda^i/i^* \geq H(i)\lambda^i$. Therefore, we have

$$\frac{\lambda}{i^*} \frac{\partial F(q, C, \vec{\lambda})}{\partial \lambda} = \frac{\lambda}{i^*} \sum_{i=1}^{2C-2} iH(i)\lambda^{i-1} = \sum_{i=1}^{i^*-1} \frac{i}{i^*} H(i)\lambda^i + \sum_{i=i^*}^{2C-2} \frac{i}{i^*} H(i)\lambda^i$$

$$\geq \sum_{i=1}^{2C-2} H(i)\lambda^i = F(q, C, \vec{\lambda}) - H(0)$$

$$> F(q, C, \vec{\lambda}) > 0,$$

where the second-to-last inequality follows because $H(0) < 0$. Thus, whenever $F(q, C, \vec{\lambda}) > 0$, we also have $\partial F(q, C, \vec{\lambda})/\partial \lambda > 0$, which (as argued above) implies that $\mathcal{A}$ is monotone in $\lambda$.

Thus, it only remains to prove the claim. Note first that, for $0 \leq i < C - 1$, $H(i) < 0$, since $\lambda^i$ appears only in the $-4q\Lambda_{C-1}^2$ term. Now, if $F(q, C, \vec{\lambda}) > 0$, then it cannot be that $H(i) \leq 0$ for all $0 \leq i \leq 2C - 2$ and, therefore, there must exist a minimal $C - 1 \leq i^* \leq 2C - 2$ such that $H(i^*) > 0$ and, for all $i < i^*$, $H(i) \leq 0$. Thus, to prove the claim, it suffices to show that $H(i) > 0$ for all $i \geq i^*$.

Now, for $C - 1 \leq i \leq 2C - 2$, observe that the coefficient $H(i)$ for $\lambda^i$ is given by

$$H(i) = (1+q)^2 \left[ \frac{1}{(C-1)!(i-C+1)!} - \frac{1}{C!(i-C)!} \right] - 4q \sum_{k=i-C+1}^{C-1} \frac{1}{k!} \frac{1}{(i-k)!}$$

$$= \frac{(1+q)^2 - ((C!(i-C+1)!)/(2C-i-1))4q\sum_{k=i-C+1}^{C-1}(1/k!)1/(i-k)!}{(C!(i-C+1)!)/(2C-i-1)}.$$

It is easy to see that then $H(i) > 0$ if and only if

$$q \geq G(i) \doteq \frac{4C!(i-C+1)!}{2C-i-1} \sum_{k=i-C+1}^{C-1} \frac{1}{k!}\frac{1}{(i-k)!} - 2.$$

Now, note that if

(4.6) $$G(i+1) \leq G(i) \qquad \text{for } i = C-1, \ldots, 2C-3,$$

then $H(i^*) > 0$, which implies $q \geq G(i^*)$, which, in turn, implies that, for all $i > i^*$, $q \geq G(i)$ and, therefore, $H(i) > 0$. Thus, in order to establish the claim, it suffices to prove (4.6).

To establish (4.6), we will find it convenient to first change variables, namely, to set $j \doteq 2C - 1 - i$ and $m \doteq C - k$ for $k = 1, \ldots, C-1$. Then, for $j = 1, \ldots, C$,

$$G(j) = \frac{4C!(C-j)!}{j} \sum_{k=C-j}^{C-1} \frac{1}{k!}\frac{1}{(2C-j-1-k)!} - 2$$

$$= \frac{4C!(C-j)!}{j} \sum_{m=1}^{j} \frac{1}{(C-m)!}\frac{1}{(C-j-1+m)!} - 2$$

$$= \frac{4C}{j} \sum_{m=1}^{j} \frac{(C-1)!(C-j)!}{(C-m)!(C-j-1+m)!} - 2$$

$$= \frac{4C}{j} \sum_{m=1}^{j} T(j,m) - 2,$$

where we define

$$T(j,m) \doteq \frac{(C-1)!(C-j)!}{(C-m)!(C-j-1+m)!}.$$

We now need to show that $G(j) \leq G(j+1)$ for $j = 1, \ldots, C-2$.

We first summarize some properties of $T(j,m)$ that will be used in the sequel. It is easy to see that $T(j,m) = T(j, j-m+1)$ and, after some manipulation, that $T(j, m+1) = (C-m)T(j,m)/(C-j+m)$. Thus, $T(j, m+1) > T(j,m)$ if and only if $(C-m)/(C-j+m) > 1$, that is, if and only if $m < j/2$. Similarly, $T(j, m+1) < T(j,m)$ for $m > j/2$, and $T(j, m+1) = T(j,m)$ for $m = j/2$. Consequently, for a given value of $j$, $T(j,m)$ has its maximum value at $m = j/2$ when $j$ is even, in which case it takes the value $(C-1)!(C-j)!/[(C-j/2)!(C-j/2-1)!]$, whereas, when $j$ is odd, the maximum $(C-1)!(C-j)!/[(C-(j+1)/2)!(C-(j+1)/2)!]$ is achieved at $m = (j + 1)/2$. It is also straightforward to check that $T(j + 1, m) = (C - j - 1 + m)T(j, m)/(C - j) \geq T(j, m)$, since $m$ is positive.

Now, let us consider $G(j)$ again. If $j$ is even, we split the sum as follows:

$$G(j) = \frac{4C}{j}\left[\sum_{m=1}^{j/2} T(j, m) + \sum_{m=1}^{j/2} T(j, j - m + 1)\right] - 2,$$

where we have used the fact that $T(j, m) = T(j, j - m + 1)$. Then, again making repeated use of the properties of $T(j, m)$ listed above, note that

$$G(j + 1) = \frac{4C}{j + 1}\sum_{m=1}^{j+1} T(j + 1, m) - 2$$

$$= \frac{4C}{j + 1}\left[\sum_{m=1}^{j/2} T(j + 1, m) + T(j + 1, j/2 + 1)\right.$$

$$\left. + \sum_{m=1}^{j/2} T(j + 1, j - m + 2)\right] - 2$$

$$\geq \frac{4C}{j + 1}\left[\sum_{m=1}^{j/2} T(j, m) + T(j + 1, j/2 + 1)\right.$$

$$\left. + \sum_{m=1}^{j/2} T(j, j - m + 1)\right] - 2$$

$$= \frac{4C}{j + 1}\left[\sum_{m=1}^{j} T(j, m) + T(j + 1, j/2 + 1)\right] - 2.$$

Now,

$$T(j + 1, j/2 + 1) = \frac{(C - 1)!(C - j - 1)!}{(C - j/2 - 1)!(C - j/2 - 1)!}$$

$$= \frac{C - j/2}{C - j}T(j, j/2)$$

$$> T(j, j/2).$$

On the other hand, $T(j, j/2)$ is the maximal value of $T(j, m)$, and so

$$T(j + 1, j/2 + 1) \geq T(j, m)$$

$$\implies \quad jT(j + 1, j/2 + 1) \geq \sum_{m=1}^{j} T(j, m)$$

$$\implies \quad T(j + 1, j/2 + 1) \geq \frac{1}{j}\sum_{m=1}^{j} T(j, m).$$

Hence,

$$G(j+1) \geq \frac{4C}{j+1} \left[ \sum_{m=1}^{j} T(j,m) + \frac{1}{j} \sum_{m=1}^{j} T(j,m) \right] - 2$$

$$= \frac{4C}{j} \sum_{m=1}^{j} T(j,m) - 2$$

$$= G(j).$$

If $j$ is odd, a similar argument gives the result. We omit the details. Thus, for $j = 1, \ldots, C-1$, we have shown that $G(j+1) \geq G(j)$ with equality if and only if $j = 1$. Returning to our original variable $i = 2C - 1 - j$, we therefore have $G(i+1) < G(i)$ for $C \leq i \leq 2C - 2$ and $G(2C-2) = G(2C-3)$, which completes the proof of (4.6) and, therefore, of the claim and theorem. $\square$

REMARK.    It is worthwhile to note that the set $\mathcal{A}$ is *not* monotone in $\lambda$ (in the same sense as specified in Theorem 4.3) when the weight vector $\vec{\lambda}$ is given by (1.5) for some $\lambda > 0$. To see this, consider, for instance, the case $C = 2$. Then, Condition A reduces to the inequality $(1+q)^2 \lambda - 4q(1+\lambda)^2 > 0$, which clearly fails to hold for $\lambda$ sufficiently large. In particular, when $q = 14$ (with $C = 2$), Condition A holds only for the weight vectors $\vec{\lambda}$ associated with $\lambda \in (49/56, 64/56)$. In other words, for $q = 14$, $C = 2$, there are multiple Gibbs measures for the parameters $(q, C, \vec{\lambda}, \nu)$ if and only if $\lambda \in (49/56, 64/56)$ and $\nu$ lies in a nonempty open interval $(\nu_-(\lambda), \nu_+(\lambda))$, while, for all other values, there is a unique Gibbs measure.

## 5. Numerical examples.
In this section, we illustrate our results with some examples.

An example of nonmonotonicity of the phase transition with respect to the multicast arrival rate $\nu$ can be seen in Figure 1, which plots the blocking probability for multicast calls when $C_E = C = 2$, $C_V = 1$, $\lambda = 0.75$ and $q = 10$, with weights $\lambda_i = \lambda^i / i!$. It can be clearly seen that there is an interval $(\nu_-, \nu_+)$ within which the fixed point is unstable, while outside this interval it is stable. This should be contrasted with the case where $C = 1$ and $C_V = 1$, where it is known [26] that the phase transition is monotone in $\nu$. In this example, we see that, as a result of the control, a difference in blocking probabilities in the even and odd sublattices of the tree, though it exists, is not significant and, hence, the unfairness due to heterogeneous blocking probabilities is not severe.

Figure 2 plots $\lambda$ against $\nu_-$ and $\nu_+$ when $C_E = C = 2$, $C_V = 1$ and $q = 6$. Note that, for these parameter values, there is no phase transition when $\lambda < 6$. The region enclosed within the curve is the parameter region where the Gibbs measure is not unique. One very interesting feature that is not immediately apparent from this plot is that, for $\lambda \in (6.000, 6.103)$, $\nu_-$ decreases as $\lambda$ increases. In other words, the point at which multiple Gibbs measures appear may *decrease* as $\lambda$ increases. This contrasts with the case $C = 1$, where an increasing $\lambda$ always increases the point at which multiple Gibbs measures appear.

Figure 3 is a similar plot for $C_E = C = 2$ and $C_V = 1$, but this time with $q = 14$ and weights $\lambda_i = \lambda^i$. This illustrates the result that, for these weights, the occurrence of multiple Gibbs measures is *not* monotone in $\lambda$; multiple Gibbs measures only occur for $\lambda \in (49/56, 64/56)$.

The phenomenon of nonmonotonicity appears not to be restricted to the controlled model, with $C_V < C$. We conclude this section with the plot in Figure 4 of the multicast blocking probabilities when $C_E = C_V = C = 2$, $q = 10$ and $\lambda = 0.72$ for $\nu \in (10, 60)$, which illustrates this point.

## 6. Concluding remarks.

In this paper, we analyzed an idealized model of multicasting on a regular tree network with a simple control that limits the number of multicast calls that can be centered at any node. With this control in place, we showed that a phase transition may exist in the infinite network and that it is nonmonotone in $\nu$, the arrival rate of the multicast calls.

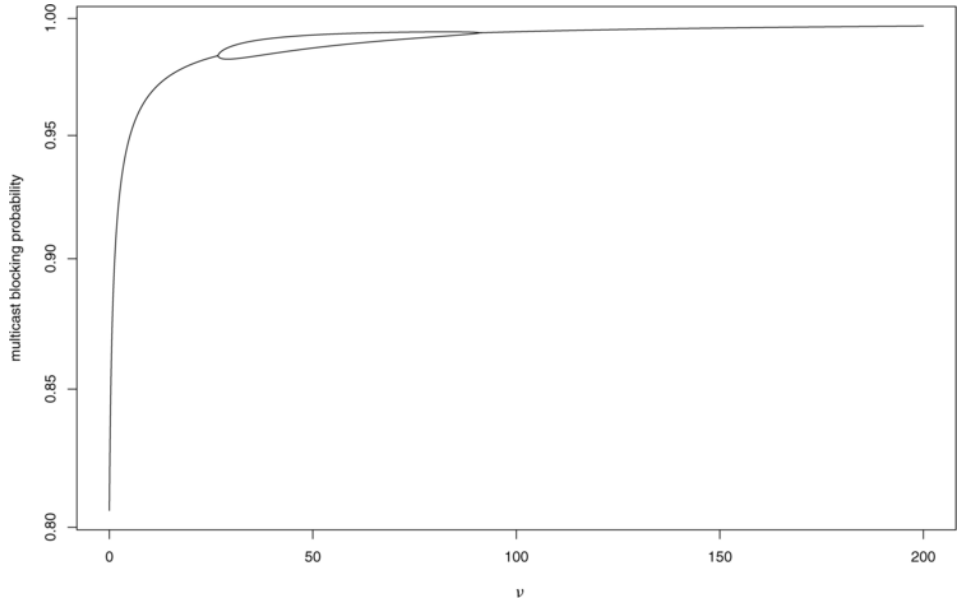

FIG. 1. *Multicast blocking probability for weights* $\lambda_i = \lambda^i/i!$ *($C_E = C = 2$, $q = 10$, $C_V = 1$, $\lambda = 0.75$).*

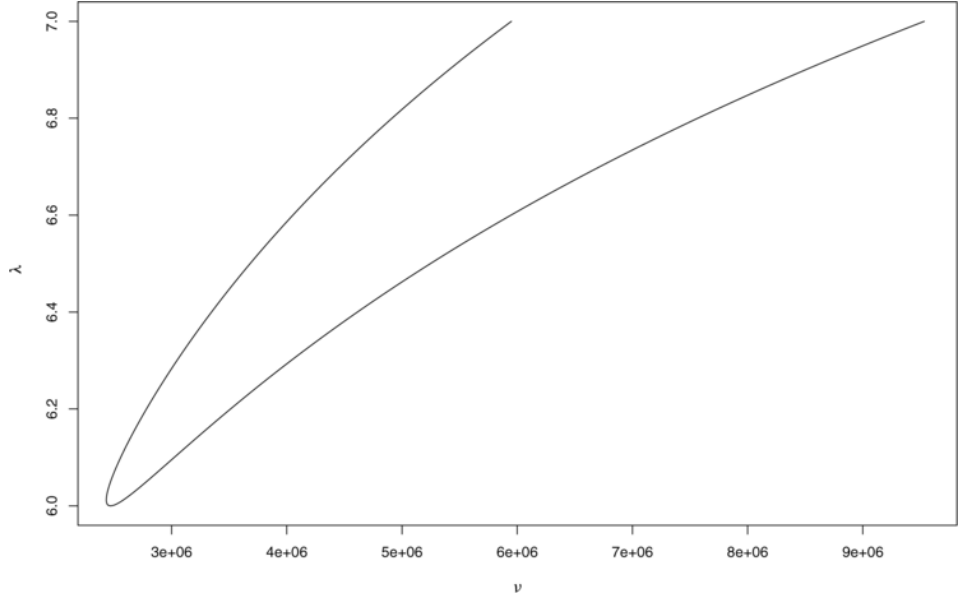

Fig. 2. *The phase transition region for weights $\lambda_i = \lambda^i/i!$ ($C_E = C = 2$, $C_V = 1$, $q = 6$).*

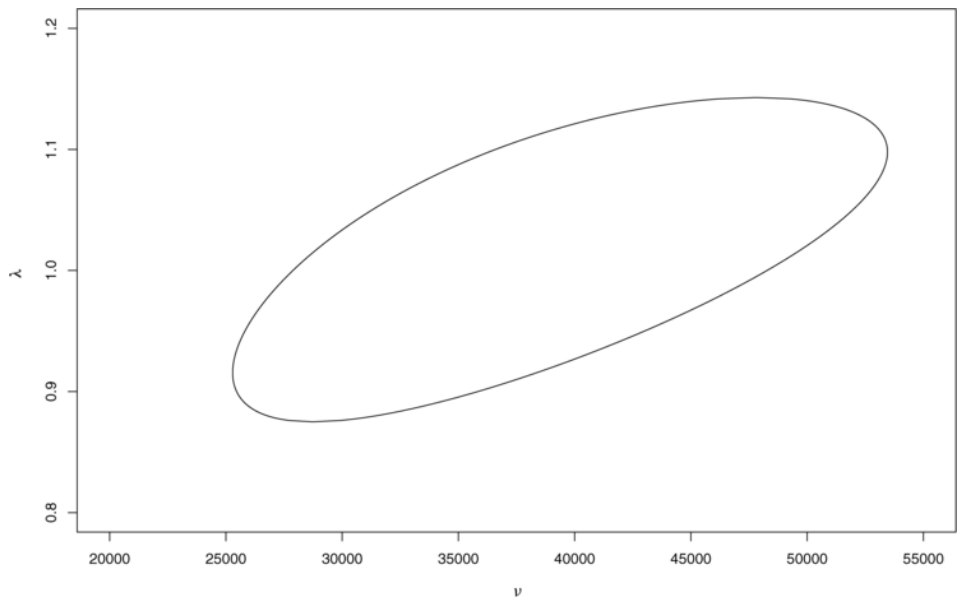

Fig. 3. *The phase transition region. $C = 2$, $q = 14$, $C_V = 1$, $\lambda_i = \lambda^i$.*

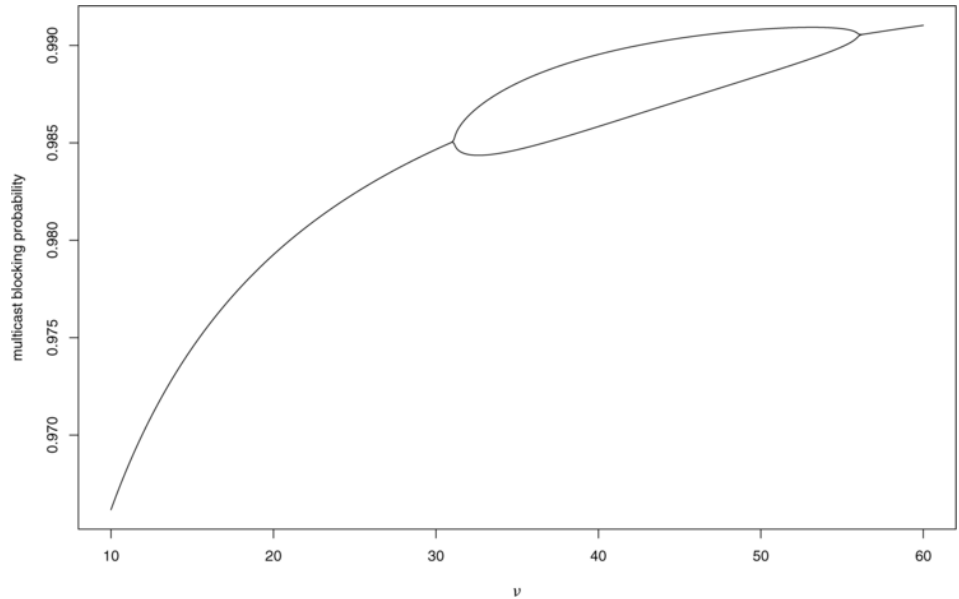

Fig. 4. *Multicast blocking probability for a model without controls.* $C_V = C_E = C = 2$, $q = 10$, $\lambda_i = \lambda^i/i!$, $\lambda = 0.72$.

It is in general a nontrivial task to identify when a phase transition is monotone with respect to the weights. For example, the identification of graphs for which the phase transition of even the hard-core model is monotone remains an open problem and, in particular, this question has not been resolved even for the $d$-dimensional lattice (see [16] for some partial results in this direction). Nonmonotonicity of phase transitions for a hard core model on a tree network with added hanging links was observed by Brightwell, Haggstrom and Winkler [4]. In this work, we presented an example of nonmonotonicity with respect to a parameter that has arisen naturally in the context of controls on multicast calls in a regular communications network.

In order to gain insight into the nature of the phase transitions studied here, it is useful to examine the different ways in which the two types of calls (i.e., unicast and multicast calls) can be packed into the network. A nonrigorous, but intuitive, explanation of the nonmonotonicity of the phase transition with respect to the multicast arrival rates is as follows. For simplicity we concentrate on the case when $C_E = C = 2$ (with $C_V = 1$). Then, for all sufficiently low $\nu$, the network behaves more or less as it would in the absence of multicast calls, in which case (since unicast calls at different links do not interact with each other) there is clearly a unique Gibbs measure. This unique measure can be thought to roughly correspond to a packing of two unicast calls throughout the network, which is clearly homogeneous.

However, as $\nu$ increases, there is increasing demand for capacity by multicast calls, leading to competition between multicast calls and unicast calls. As a consequence, the unicast calls are now restricted to occupy only one unit of capacity throughout the network, resulting in the multicast calls effectively experiencing a tree network with capacity 1. In this case, as is well known [18], there are two phases, with one corresponding to a packing of one multicast call on each node on the even subtree, and the other to a packing of each multicast call on the odd subtree. As the multicast arrival rate increases yet further (relative to the unicast rate), the unicast calls are unable to capture even one unit of capacity throughout the network and, thus, one obtains once again a single homogeneous Gibbs measure corresponding to a packing of a multicast call at each node of the network. Similar reasoning leads us to the following conjecture, which is supported by numerical results (see Figure 4 of this paper and Section 5 of [21]).

CONJECTURE. For the *uncontrolled* unicast-multicast model with $C_E = C_V = C = 2$, there exist phase transitions that are nonmonotone in the multicast arrival rate $\nu$.

Recent years have witnessed an explosion of research papers devoted to the analysis of phase transitions in statistical physics-type models in various fields, with the majority arising from problems in combinatorics and the theory of computing (see, e.g., [3, 6] and references therein). This work, along with [13, 17, 24] and [26], demonstrates that Gibbs measures associated with loss networks also exhibit several interesting phenomena that are worthy of further study.

**Acknowledgments.** Kavita Ramanan would like to thank the NZIMA for funding a visit to Auckland, New Zealand, the Statistics Department at the University of Auckland for its hospitality during this visit, and also the IMA in Minnesota for its funding and hospitality. Ilze Ziedins would like to thank the NZIMA, the IMA in Minnesota and the Mittag-Leffler Institute in Stockholm for their funding and the latter two for their hospitality. This research also benefited from visits to Bell Labs, as well as discussions with Stan Zachary. We also wish to thank the referee for a very careful reading.

B. LUEN
DEPARTMENT OF STATISTICS
UNIVERSITY OF CALIFORNIA, BERKELEY
BERKELEY, CALIFORNIA 94720
USA
E-MAIL: bradluen@stat.berkeley.edu

K. RAMANAN
DEPARTMENT OF MATHEMATICAL SCIENCES
CARNEGIE MELLON UNIVERSITY
PITTSBURGH, PENNSYLVANIA 15213
USA
E-MAIL: kramanan@math.cmu.edu

I. ZIEDINS
DEPARTMENT OF STATISTICS
UNIVERSITY OF AUCKLAND
PRIVATE BAG 92019
AUCKLAND
NEW ZEALAND
E-MAIL: ilze@stat.auckland.ac.nz